\documentclass[a4paper,11pt]{article}
\usepackage[latin1]{inputenc}
\usepackage{amsfonts}
\usepackage{amsthm}
\usepackage{amsmath}
\usepackage{amsfonts}
\usepackage{latexsym}
\usepackage{amssymb}
\usepackage{dsfont}

\newtheorem{fed}{Definition}[section]
\newtheorem{teo}[fed]{Theorem}
\newtheorem*{teo*}{Theorem}
\newtheorem{lem}[fed]{Lemma}

\newtheorem{pro}[fed]{Proposition}
\theoremstyle{definition}
\newtheorem{rem}[fed]{Remark}

\newtheorem{num}[fed]{}
\def\ID{Let $(\cF_0 \coma \ca)$ be initial data
for the CP} 
\def\ga{\gamma}

\def\n0{n_{ \text{\rm \tiny o}}}

\newcommand{\IN}[1]{\mathbb {I} _{#1}}
\def\In{\mathbb {I} _n}
\def\IM{\mathbb {I} _m}
\def\suml{\sum\limits}
\def\QEDP{\tag*{\QED}}

\def\bce{\begin{center}}
\def\ece{\end{center}}
\newcommand{\trivial}{\{0\}}
\DeclareMathOperator{\FP}{FP\,}

\def\ds{\displaystyle} 
\def\subim{_{i\in \IN{n}}\,}
\def\py{\peso{and}}
\def\rk{\text{\rm rk}}
\def\noi{\noindent}
\def\cF{\mathcal F}
\def\cG{\mathcal G}
\def\QED{\hfill $\square$}
\def\EOE{\hfill $\triangle$}

\newcommand{\peso}[1]{ \quad \text{ #1 } \quad }
\def\uno{\mathds{1}}
\def\bm{\left[\begin{array}}
\def\em{\end{array}\right]}
\def\ben{\begin{enumerate}}
\def\een{\end{enumerate}}
\def\bit{\begin{itemize}}
\def\eit{\end{itemize}}
\def\barr{\begin{array}}
\def\earr{\end{array}}
\def\igdef{\ \stackrel{\mbox{\tiny{def}}}{=}\ }

\def\k{n}

\def\la{\lambda}

\def\N{\mathbb{N}}
\def\R{\mathbb{R}}
\def\C{\mathbb{C}}

\def\cC{\mathcal{C}}

\def\cH{\mathcal{H}}
\def\cK{\mathcal{K}}

\def\cS{{\cal S}}

\def\cB{{\cal B}}

\def\cV{{\cal F}}

\def\cW{{\cal G}}
\def\ca{\mathbf{a}}
\def\cb{\mathbf{a}}

\def\vacio{\varnothing}
\def\orto{^\perp}
\def\inc{\subseteq}

\def\sii{ if and only if }
\def\inv{^{-1}}

\def\api{\langle}
\def\cpi{\rangle}

\def\ua{^\uparrow}
\def\da{^\downarrow}

\DeclareMathOperator{\Preal}{Re} 
 \DeclareMathOperator{\tr}{tr}
\DeclareMathOperator{\gen}{span}
\DeclareMathOperator{\convf}{Conv (\R_{\ge0})}
\DeclareMathOperator{\convfs}{Conv_s (\R_{\ge0})}
\DeclareMathOperator{\leqp}{\leqslant}

\def\RS{\mathbf{F} }

\def\RSV{\cF= \{f_i\}_{i\in \, \IN{n}}}

\def\coma{\, , \, }
\def\elnu{\nu_f(\la\coma \ca)}
\def\elmu{\mu_f(\la\coma \ca)}
\def\nuel{\nu(\la\coma \ca)}
\def\muel{\mu(\la\coma \ca)}

\newcommand{\hil}{\mathcal{H}}
\newcommand{\op}{L(\mathcal{H})}
\newcommand{\lhk}{L(\mathcal{H} \coma \mathcal{K})}

\newcommand{\posop}{L(\mathcal{H})^+}

\def\H{{\cal H}}
\def\glh{\mathcal{G}\textit{l}\,(\cH)}

\newcommand{\mat}{\mathcal{M}_d(\mathbb{C})}

\newcommand{\matsad}{\mathcal{H}(d)}

\newcommand{\matud}{\mathcal{U}(d)}

\newcommand{\matpos}{\mat^+}

\newcommand{\matinvd}{\mathcal{G}\textit{l}\,(d)}

\newcommand{\matrec}[1]{\mathcal{M}_{#1} (\mathbb{C})}

\def\beq{\begin{equation}}
\def\eeq{\end{equation}}

\def\pausa{\medskip\noi}

\begin{document}

\title{ {\bf Optimal frame completions with prescribed norms for majorization
}}
\author{Pedro G. Massey, Mariano A. Ruiz  and Demetrio Stojanoff\thanks{Partially supported by CONICET 
(PIP 5272/05) and  Universidad Nacional de La PLata (UNLP 11 X472).} }
\author{P. G. Massey, M. A. Ruiz and D. Stojanoff \\ {\small Depto. de Matem\'atica, FCE-UNLP,  La Plata, Argentina
and IAM-CONICET  \footnote{e-mail addresses: massey@mate.unlp.edu.ar , mruiz@mate.unlp.edu.ar , demetrio@mate.unlp.edu.ar}
}}
\date{}
\maketitle

\begin{abstract}  Given a finite sequence of vectors $\mathcal F_0$ in $\C^d$ we characterize in a complete and explicit way the optimal completions of $\mathcal F_0$ obtained by adding a finite sequence of vectors with prescribed norms, where optimality is measured with respect to majorization (of the eigenvalues of the frame operators of the completed sequence). Indeed, we construct (in terms of a fast algorithm) a vector - that depends on the eigenvalues of the frame operator of the initial sequence $\cF_0$ and the sequence of prescribed norms - that is a minimum for majorization among all eigenvalues of frame operators of completions with prescribed norms. Then, using the eigenspaces of the frame operator of the initial sequence $\cF_0$ we describe the frame operators of all optimal completions for majorization. Hence, the concrete optimal completions with prescribed norms can be obtained using recent algorithmic constructions related with the Schur-Horn theorem.

The well known relation between majorization and tracial inequalities with respect to convex functions allow to describe our results in the following equivalent way: given a finite sequence of vectors $\mathcal F_0$ in $\C^d$ we show that the completions with prescribed norms that minimize the convex potential induced by a strictly convex function are structural minimizers, in the sense that they do not depend on the particular choice of the convex potential.     

\end{abstract}



\date{}

\tableofcontents

\section{Introduction} 

A finite sequence of vectors $\cF=\{f_i\}_{i\in \IN{n}}$ in $\C^d$ is a frame for $\C^d$ if the sequence spans $\C^d$. It is well known that finite frames provide (stable) linear encoding-decoding schemes. As opposed to bases, frame are not subject to linear independence; indeed, it turns out that the redundancy allowed in finite frames can be turned into robustness of the transmission scheme that they induce, which make frames a useful device for transmission of signals through noisy channels (see \cite{Bod,Pau,BodPau,caskov,HolPau,LoHanagre,LeHanagre}). 
 
 On the other hand, the so-called tight frames allow for redundant linear representations of vectors  that are formally analogous to the linear representations given by orthonormal basis; this feature makes tight frames a distinguished class of frames that is of interest for applications. 
 In several applications we would like to consider tight frames that have some other prescribed properties leading to what is known in the literature as  frame design problems \cite{Illi,CFMP,Casagregado,CMTL,ID,DFKLOW,FWW,KLagregado}. It turns out that in some cases it is not possible to find a frame fulfilling the previous demands.
 
An alternative approach to deal with the construction of frames with prescribed parameters and nice associated reconstruction formulas was posed in \cite{BF} by Benedetto and Fickus; they defined a functional, called the frame potential, and showed that minimizers of the frame potential (within a convenient set of frames) are the natural substitutes of tight frames with prescribed parameters (see also \cite{Phys,FJKO,JOk,MR} and \cite{casafick3,MRS,MRS2} for related problems in the context of fusion frames). Moreover, in \cite{MR} it is shown that minimizers of the frame potential under suitable restrictions (considered in the literature) are structural minimizers in the sense that they coincide with minimizers of more general {\it convex potentials} (see Section \ref{basic}).

Recently, the following frame completion problem was posed in \cite{FMP} (in the vein of \cite{BF}): given
an initial sequence $\cF_0$ in $\C^d$ and a sequence of positive numbers $\cb$ then compute the sequences $\cG $ in $\C^d$ whose elements have norms given by the sequence $\cb$ and such that 
the completed sequence $\cF=(\cF_0 \coma \cG )$ minimizes the so-called mean square error (MSE) of $\cF$, which is a (convex) functional (see also \cite{CCHK, FWW,MR0} for completion problems for  frames). In this setting, the initial sequence of vectors can be considered as a checking device for the measurement, and therefore we search for a complementary set of measurements (given by vectors with prescribed norms) in such a way that the complete set of measurements is optimal with respect to the MSE. Notice there are other possible (convex) functionals that we could choose to minimize such as, for example, the frame potential. Therefore,
 a natural extension of the previous problem is: given a  functional defined on the set of frames, compute the frame completions with prescribed norms that minimize this functional. Moreover, this last problem raises the question of whether the completions that minimize these functionals coincide i.e., whether the minimizers are structural in this setting.

A first step towards the solution of the general version 
of the completion problem was made in \cite{MRS4}. There we showed 
that under certain hypothesis (feasible cases, see 
Section \ref{TFC}), optimal frame completions 
with prescribed norms are structural (do not depend on the particular choice of the convex 
functional), as long as we consider 
convex potentials, that contain the MSE and the frame potential. 
On the other hand, it is easy to show examples in which the previous 
result does not apply (non-feasible cases); in these cases the optimal 
frame completions with prescribed norms were not known even for the 
MSE nor the frame potential. Recently, in some feasible cases  
the set of all optimal frame completions is characterized in \cite {Pot, FMPS}.

In \cite{MRS5} we considered the structure of completions that minimize a fixed convex potential 
(non feasible case). There, we showed that the eigenvalues of optimal completions with respect to a fixed convex potential are uniquely determined by the solution of an optimization problem in a compact convex subset of $\R^d$ for a convex objective function that is associated to the convex potential in a natural way.
Then, we showed an important geometrical feature of optimal completions $\cF=(\cF_0,\cG)$ for a fixed convex potential, namely that the vectors in the completion $\cG$ are eigenvectors of the frame operator of the completed sequence $\cF$ (see Section \ref{sec3} for a detailed exposition of these results). Based on these facts, we developed an algorithm that allowed us to compute the solutions of the completion problem for small dimensions. In this setting we conjectured some properties of the optimal frame completions in the general case, based on common features of the solutions of several examples obtained by this algorithm (see Section \ref{sec The main result} for a detailed description of these conjectures).

In this paper, building on our previous work \cite{MRS4} and \cite{MRS5}, we give a complete and explicit description of the spectral and geometrical structure of optimal completions with prescribed norms with respect to a convex potential induced by a strictly convex function. Our approach is constructive and allows to develop a fast and effective algorithm that computes the spectral structure of optimal completions. As we shall see, 
given an initial sequence $\cF_0$ in $\C^d$ and a sequence of positive numbers $\cb$, 
both the spectral and geometrical structure of optimal completions depend only on the frame operator of $\cF_0$ and $\cb$, but they do not depend on the particular choice of the convex potential. Hence, we show that in the general case the minimizers of convex potentials (induced by strictly convex functions) are structural. 

In order to obtain the previous results, we begin by proving the properties of general optimal completions conjectured in \cite{MRS5}. These properties (that are structural, in the sense that they do not depend on the convex potential) are then used to compute several other structural parameters - that involve the notion of feasibility developed in \cite{MRS4} - that completely describe the spectral structure of optimal completions. As a consequence of this description, we conclude that optimal solutions have the same eigenvalues and hence, the eigenvalues of optimal completions are minimum for the so-called majorization preorder. Moreover, all the parameters 
involved in the description of the spectral structure of optimal completions
can be computed in terms of fast algorithms.
With the spectral data and results from \cite{MRS4} we completely describe the set positive matrices that correspond to the frame operators of sequences $\cG$ with norms prescribed by $\cb$ and such that $\cF=(\cF_0,\cG)$ are optimal. Finally, some optimal 
completions $\cG$ can be also 
effectively computed by using recent results from
\cite{CFMP} (see also \cite{ID} and \cite {FMPS}) that implements the Schur-Horn theorem.


The paper is organized as follows. In Section \ref{Pre} we briefly recall the basic framework of finite frame theory, the notion of submajotization 
- that will play a central role in this note - and the relation of submajorization with tracial inequalities involving convex functions. 
In section \ref{problemon} we describe the context of our main problem - namely, optimal completions with prescribed norms, where optimality is described in terms of majorization - and give a detailed account of several related results that were developed in our previous works \cite{MRS4} and \cite{MRS5} that we shall need in the sequel, in a way suitable for this note; in particular, we include a new construction of the spectra of optimal completions in the feasible cases. In Section \ref{sec The main result} we introduce new structural parameters - that can be efficiently computed in terms of explicit algorithms - and show how to give a complete description of the spectra of optimal completions for strictly convex potentials, in terms of these parameters in the general case. This allow us to show that the spectra of such optimal completions do not depend on the choice of strictly convex potential, so that minimizers are then structural. The proofs of the technical results of this section is presented in Section 5. 
In particular, we settle in the affirmative some features of the structure of optimal completions for strictly convex potentials that were conjectured in \cite{MRS5}. As a byproduct we also settle in the affirmative a conjecture on local minimizers of strictly convex potentials with prescribed norms posed in \cite{MR}.

%
%
%

\section{Preliminaries}\label{Pre}

In this section we describe the basic notions that we shall consider throughout the paper. In Section \ref{sec 2.1} we describe some general notations and terminology. In Section \ref{basic} we describe some basic notions and facts of frame theory and we recall the notion of convex potential from \cite{MR}. In Section \ref{subsec 2.2} we describe some aspects of submajorization that we shall need in the sequel.


\subsection{General notations.}\label{sec 2.1}
Given $m \in \N$ we denote by $\IM = \{1, \dots , m\} \inc \N$ and 
$\uno = \uno_m  \in \R^m$ denotes the vector with all its entries equal to $1$. 
For a vector $x\in \R^m$ we denote by 
$\tr \, x = \sum_{i\in \IN{m}} \, x_i$ 
and by $x^\downarrow$ (resp. $x^\uparrow$) the rearrangement
of $x$ in  decreasing (resp. increasing) order. We denote by
$(\R^m)^\downarrow = \{ x\in \R^m : x = x^\downarrow\}$ the set of downwards ordered vectors, and similarly $(\R^m)\ua $. 

\pausa
Given $\cH \cong \C^d$  and $\cK \cong \C^n$, we denote by $\lhk $ 
the space of linear transformations $T : \cH \to \cK$. 
If $\cK = \cH$ we denote by $\op = L(\cH \coma \cH)$, 
by $\glh$ the group of all invertible operators in $\op$, 
 by $\posop $ the cone of positive operators and by
$\glh^+ = \glh \cap \posop$. 
If $T\in \op$, we  denote by   
$\sigma (T)$ the spectrum of $T$, by $\rk\, T= \dim R(T) $  the rank of $T$,
and by $\tr T$ the trace of $T$.

\pausa
If $W\inc \cH$ is a subspace we denote by $P_W \in \posop$ the orthogonal 
projection onto $W$. 
Given $x\coma y \in \cH$ we denote by $x\otimes y \in \op$ the rank one 
operator given by 
$x\otimes y \, (z) = \api z\coma y\cpi \, x$ for every $z\in \cH$. Note that
if $\|x\|=1$ then $x\otimes x = P_{\gen\{x\}}\,$.

\pausa 
By fixing orthonormal basis's (ONB's) 
of the Hilbert spaces involved, we shall identify operators with 
matrices, using the following notations: 
by $\matrec{n,d} \cong L(\C^d \coma \C^n)$ we denote the space of complex $n\times d$ matrices. 
If $n=d$ we write $\mat = \matrec{d,d}$ ;   
$\matsad$ is the $\R$-subspace of selfadjoint matrices,  
$\matinvd$ the group of all invertible elements of $\mat$, $\matud$ the group 
of unitary matrices in $\mat$, 
$\matpos$ the cone of positive semidefinite
matrices, and $\matinvd^+ = \matpos \cap \matinvd$.

\pausa
Given $S\in \matpos$, we write $\la(S) = \la\da(S)\in 
(\R_{\geq 0}^d)^\downarrow$ the 
vector of eigenvalues of $S$ - counting multiplicities - arranged in decreasing order. 
Similarly we denote by $\la\ua(S) \in (\R_{\geq 0}^d)\ua$ the reverse ordered 
vector of eigenvalues of $S$. 
If $ \la= (\la_i)_{i\in\IN{d}}\in\R_{\geq 0}^d$ (not necessarily ordered), 
 a system $\cB=\{h_i\}_{i\in \IN{d}} \inc \C^d$ is a 
``ONB of eigenvectors for $S\coma \la \,$" if it is an  
orthonormal basis for $\C^d$ such that 
$S\,h_i=\lambda_i\,h_i$ for every $i\in \IN{d}\,$. 
In other words, an  orthonormal basis 
\beq\label{BON S la}
\mbox{$\cB=\{h_i\}_{i\in \IN{d}} $ \ \ \ is a 
``ONB of eigenvectors for $S\coma \la \,$"} \iff 
S = \sum_{i\in \IN{d}} \, \la_i \cdot h_i \otimes h_i \ .
\eeq

\subsection{Basic framework of finite frames }\label{basic}

In what follows we consider $(\k,d)$-frames. See 
\cite{BF,TaF,Chr,HLmem,MR} for detailed expositions of several aspects of this notion. 

\pausa
Let $d, n \in \N$, with $d\le n$. Fix a Hilbert space $\hil\cong \C^d$. 
A family $ \RSV \in  \cH^n $  is an 
$(\k,d)$-frame for $\cH$  if there exist constants $A,B>0$ such that
\beq\label{frame defi} A\|x\|^2\leq \sum_{i=1}^n |\left \langle x \, , f_i\right \rangle|^2\leq B \|x\|^2 \peso{for every} x\in \hil \ .
\eeq
The {\bf frame bounds}, denoted by  $A_\cF, B_\cF$ are the optimal constants in \eqref{frame defi}. If $A_\cF=B_\cF$ we call $\cF$ a tight frame.
Since $\dim \hil<\infty$, a family  $\RSV$  is an 
$(\k,d)$-frame 
 \sii $\gen\{f_i: i \in \In \} = \cH$.  
We shall denote by $\RS = \RS(n \coma d)$ the set of all $(\k,d)$-frames for $\cH$. 

\pausa
Given $\RSV \in  \cH^n $, the operator $T_\cV \in L(\hil\coma\C^n)$ defined
by 
 \beq \ T_\cV\, x= \big( \,\api x \coma f_i\cpi\,\big)  \subim 
\ , \peso{for every} x\in \cH \,
\eeq
is the {\bf analysis} operator of $\cF$.  Its adjoint $T_\cV^* \in L(\C^n\coma \cH)$  is called the {\bf synthesis} operator and is given by 
 $T_\cV ^* \, v =\sum_{i\in \, \IN{n}} v_i\, f_i$ 
for every $v = (v_i)_{i\in\In}\in \C^n$. 
The {\bf frame operator} of $\cV$ is 
$$
\barr{rl}
S_\cV &= T_\cV^*\  T_\cV = \sum_{i \in \In} f_i \otimes f_i 
\in \posop\ . \earr
$$
Notice that, if 
$\RSV \in  \cH^n $ then  
$\api S_\cV \, x\coma x\cpi \, = \sum\subim \, 
 \big|\, \api x \coma f_i\cpi \, \big|^2$ for every 
 $x\in \cH$. Hence, $\cF\in \RS(n \coma d)$ if and only if 
$S_\cF\in \glh^+$ and in this case
$A_\cV \, \|x\|^2 \, \le \, \api S_\cV \, x\coma x\cpi 
 \,  \le \, B_\cV \, \|x\|^2$ for every $ x\in \cH $.   
In particular, $A_\cV   =\la_{\min} (S_\cV) = \|S_\cV\inv \| \inv$ and $ 
\la_{\max} (S_\cV) = \|S_\cV \| = B_\cV \,$.
Moreover, $\cV $ is tight if and only if $S_\cV = 
\frac{\tau}{d}  \, I_\H\,$, where $\tau = 
\tr S_\cV = \sum\subim \|f_i\|^2 \,$.

%
\pausa 
In their seminal work \cite{BF}, Benedetto and Fickus introduced a functional defined, the so-called frame potential, given by 
$$ 
\barr{rl}
\FP(\{f_i\}_{i\in \In} )& 
=\sum_{i,\,j\,\in \In}|\api f_i\coma f_j \cpi |\,^2\ .
\earr
$$ 
One of their major results shows that tight unit norm frames - which form an important class of frames because of their simple reconstruction formulas - can be characterized as (local) minimizers of this functional among unit norm frames. Since then, there has been interest in (local) minimizers of the frame potential within certain classes of frames, since such minimizers can be considered as natural substitutes of tight frames (see \cite{Phys,MR,MRS}). Notice that, given $\cF=\{f_i\}_{i\in \IN{n}}\in  \cH^n$ then $\FP(\cF)=\tr \, S_\cF^2 
=\sum_{i\in \IN{d}}\lambda_i(S_\cF)^2$.
These remarks have motivated 
the definition of general convex potentials as follows:

\begin{fed}\label{pot generales}\rm
Let us denote by 
$$
\convf = \{ 
f:[0 \coma \infty)\rightarrow [0 \coma \infty): f   \ \mbox{ is a convex function } \ \} 
$$  
and $\convfs = \{f\in \convf : f$ is strictly convex $\}$. 
Following \cite{MR} we consider the 
(generalized) convex potential $P_f$ associated to any $f\in \convf$, given by
$$
\barr{rl}
P_f(\cF)&=\tr \, f(S_\cF) = \sum_{i\in \IN{d}}f(\lambda_i(S_\cF)\,) \peso {for} 
\cF=\{f_i\}_{i\in \IN{n}}\in  \cH^n \ , \earr
$$
where the matrix $f(S_\cF)$ is defined by means of the usual functional calculus. \EOE
\end{fed}

\pausa
As shown in \cite[Sec. 4]{MR} these convex potentials (which are related with the 
so-called entropic measures of frames) share many properties with the 
BF-frame potential. Indeed, under certain restrictions both the spectral 
and geometric structures of minimizers of these potentials coincide 
(see \cite{MR} and Remark \ref{vieja conj} below).

\begin{rem} 
The results that we shall develop in this work apply in the case of 
convex potentials $P_f$ for any $f\in \convfs$. 
Notice that this formulation does not formally include the Mean Square Error (MSE), which 
is the convex potential associated with the strictly convex function
$f:(0,\infty)\rightarrow (0,\infty)$ given by $f(x)= x\inv$, since $f$ is not defined in $0$ in this case. In order to include the MSE within our results we proceed as follows: 
we define $\tilde f:[0,\infty)\rightarrow (0,\infty]$ given by $\tilde f(x)=x^{-1}$ for $x>0$ and $\tilde f(0)=\infty$. 
Assuming that  $ x<\infty $ and $x+\infty =x\cdot \infty =\infty$ for 
every $x\in (0\coma \infty)$, it turns out that the new map $\tilde f$ is a (extended) 
strictly convex function and  all the results obtained in this paper apply to the 
convex potential induced by $\tilde f$.\EOE
\end{rem}

\subsection{Submajorization} \label{subsec 2.2}

Next we briefly describe submajorization, a notion from matrix analysis theory that will be used throughout the paper. For a detailed exposition of submajorization see \cite{Bat}.

\pausa 
 Given $x,\,y\in \R^d$ we say that $x$ is
{\bf submajorized} by $y$, and write $x\prec_w y$,  if
$$
\barr{rl}
\suml_{i=1}^k x^\downarrow _i & \le 
\suml_{i=1}^k y^\downarrow _i \peso{for every} k\in \mathbb I_d \ .\earr
$$  
If $x\prec_w y$ and $\tr x = \sum_{i=1}^dx_i=\sum_{i=1}^d y_i = \tr y$,  then we say that $x$ is
{\bf majorized} by $y$, and write $x\prec y$. In case that $x\prec y$ but $y\nprec x$ we say that $y$ majorizes $x$ strictly. If the two vectors $x$ and $y$ have different sizes, we write 
$x\prec y$ if the extended vectors (completing with zeros 
to have the same size) satisfy the previous relationship.  

\pausa
On the other hand we write 
$x \leqp y$ if $x_i \le y_i$ for every $i\in \mathbb I_d \,$.  It is a standard  exercise 
to show that $x\leqp y \implies x^\downarrow\leqp y^\downarrow  \implies x\prec_w y $. 
Majorization is usually considered because of its relation with tracial inequalities 
for convex functions. 
Indeed, if we
let $x,\,y\in \R^d$ and let $I\inc \R$ be an interval such that 
$x,\,y\in I^d$ then (see for example \cite{Bat}):
\ben 
\item $x\prec y$ $\Leftrightarrow$
$ 
\tr f(x) \igdef\suml_{i=1}^df(x_i)\leq \suml_{i=1}^df(y_i)=\tr f(y) 
$ for every convex function  $f:I\rightarrow \R$. 
\item If only $x\prec_w y$,  but the map $f:I\rightarrow \R$ is convex and increasing, then  
$\tr f(x) \le \tr f(y)$. 
\item If $x\prec y$ and $f:I\rightarrow \R$ is a strictly convex function such 
that $\tr \,f(x) =\tr \, f(y)$ then there exists a permutation $\sigma$ 
of $\IN{d}$ such that $y_i=x_{\sigma(i)}$ for $i\in \IN{d}\,$. 
\een
As a consequence of item 3. above, if $x\prec y$ strictly and $f:I\rightarrow \R$ is a strictly convex function then $\tr \,f(x) <\tr \, f(y)$: indeed, if $\tr \,f(x)=\tr \, f(y)$ then by item 3. above we would have that $y_i=x_{\sigma(i)}$, $i\in \IN{d}\,$, for a permutation $\sigma$ of $\IN{d}$ and hence that $y\prec x$.

\pausa
The notion of vector submajorization can be extended to a preorder between selfadjoint matrices as follows:
given $S_1\coma S_2\in \mathcal H(d)$ we 
say that $S_1$ is submajorized by $S_2\, $, 
and write $S_1\prec_w S_2\, $ (resp. $S_1\prec S_2\, $) if $\lambda(S_1)\prec_w \lambda(S_2)\, $ 
(resp. $\lambda(S_1)\prec \lambda(S_2)$, i.e. $S_1\prec_w S_2\, $ and $\tr \,  S_1=\tr \, S_2$). 

\begin{rem}\label{rem doublystohasctic}
Majorization between vectors in $\R^d$ is intimately related with the 
class of doubly stochastic $d\times d$ matrices, denoted by DS$(d)$. Recall that a 
$d\times d$ matrix $D \in $ DS$(d)$   if it has non-negative entries and each row sum and column sum equals 1. 
%
It is well known (see \cite{Bat}) that given $x\coma y\in \R^d$ then $x\prec y$ 
if and only if there exists $D\in$ DS$(d)$ such that $D\, y=x$. As a consequence 
of this fact we see that if $x_1\coma y_1\in \R^r$ and $x_2 \coma y_2\in \R^s$ 
are such that 
\beq\label{cachos}
x_i\prec y_i \peso{for} i=1\coma 2 
\implies x=(x_1\coma x_2)\prec y=(y_1\coma y_2)  \peso{in} 
\R^{r+s} \ .
\eeq
Indeed, if $D_1$ and $D_2$ are the doubly stochastic matrices corresponding the previous majorization relations then $D=D_1\oplus D_2\in$ DS$(r+s)$ is such that $D\, y=x$. 
\EOE
\end{rem}


\section{Optimal completions with prescribed norms}\label{problemon}

In this section we give a detailed description of the optimal completion problem 
and recall some notions and results from our previous work \cite{MRS4,MRS5}, in a way suitable for the exposition of the results herein. In particular, the exposition of the results in Section \ref{TFC} differs from that of \cite{MRS5}, since this new presentation is better suited for our present purposes.

\subsection{Presentation of the problem}
In several applied situations it is desired to construct 
a sequence $\cG$ in such a way that the  frame 
operator of $\cG$ is given by some $B\in\matpos$ and the squared norms of the frame elements 
are prescribed by a sequence of positive numbers  $\ca=(a_i)_{i\in \IN{k}}\in \R_{>0}^k\,$. 
That is, given a fixed $B\in \matpos$ and $\mathbf{a}\in \R_{>0}^k\,$, we analyze the existence (and construction) of a sequence $\cG=\{g_i\}_{i\in \IN{k}}$ such that $S_\cG=B$ 
and $\|g_i\|^2=a_i\,$, for $i\in \IN{k}\,$. This is known as the classical 
frame design problem. It  has been treated by several research groups (see for example \cite{Illi,CFMP,Casagregado,CMTL,ID,DFKLOW,FWW,KLagregado}).
In what follows we recall 
a solution of the classical frame design problem in the finite dimensional setting, in the way that it is convenient for our analysis.

\begin{pro}[\cite{Illi,MR0}]\label{frame mayo}\rm
Let $B\in \matpos$ with $\lambda(B)\in (\R_{\geq 0}^d)^\downarrow$ and let 
$\cb=(a_i)_{i\in\IN{k}} \in \R_{>0}^k\,$. Then there exists 
a sequence $\cW=\{g_i\}_{i\in\IN{k}}\in \cH^k$ with frame operator 
$S_\cW= B$ and such that $\|g_i\|^2=a_i$ for every $i\in\IN{k}\,$ 
if and only if 
$\cb\prec \lambda(B)$ (completing with zeros if $k\neq d$).
\QED
\end{pro}

\pausa
Recently, researchers have made a step forward in the classical frame design problem and have asked about the structure of {\bf optimal} frames with prescribed parameters.
Indeed, consider the following problem posed in \cite{FMP}:
let $\cH \cong \C^d$ and let $\cF_0=\{f_i\}_{i\in \IN{\n0}}
\in \cH^{\n0  }$ be a fixed (finite) sequence of vectors. 
Consider a sequence $\ca
= (a_i)_{i\in \IN{k}} \in \R_{>0}^k\,$ such
that $\rk \, S_{\cF_0} \ge d-k$ and denote by $n=\n0+k$.  Then, 
with this fixed data, the problem is to construct a sequence  
$$
\cG = \{g_i\}_{i\in \IN{k}}\in \cH^{k}  \peso{with} \|g_{i}\|^2=a_i \peso{for}
i\in \IN{k}  \ , 
$$
such that the resulting completed sequence $\cF= (\cF_0\coma \cG )\in\RS(n \coma d)$ - obtained by juxtaposition of the two finite sequences - is a frame 
 whose MSE, given by $\tr \, S_\cF^{-1}$, is minimal among 
all possible such completions.

\pausa
Note that there are other possible ways to measure robustness 
(optimality) of the completed frame $\cF$ as above. 
For example, we can consider optimal (minimizing) completions, with prescribed norms, 
for the Benedetto-Fickus' potential. In this case we search for a frame 
$\cF= (\cF_0\coma \cG )\in\RS(n \coma d)$, with
$\|g_{i}\|^2=a_i$ for $i\in \IN{k}$, and such that its frame potential 
$\FP(\cF)=\tr \,S_\cF ^2 $ is minimal  
among all possible such completions (indeed, this problem has been 
considered before in the particular case in which $\cF_0=\vacio$ 
in \cite{BF,Phys,FJKO,JOk,MR}). More generally, we can measure robustness of the completed frame $\cF=(\cF_0,\cG)$ 
in terms of general convex potentials (see Definition \ref{pot generales}).

\pausa
In order to describe the main problems we first fix the 
notation that we shall use throughout the paper.

\begin{fed}\label{cp}\rm 
Let $\cF_0=\{f_i\}_{i\in \IN{\n0  }}\in \cH^{\n0  }$ and 
$\ca= (a_i)_{i\in \IN{k}} \in (\R_{>0}^k)\da\,$ 
such that $d-\rk \, S_{\cF_0} \le k$. Define $n=\n0+k$. Then
\ben
\item In what follows we say that
 $(\cF_0 \coma \ca)$ are initial data for the completion problem (CP). 
\item 
For these data we consider the sets
$$
\cC_\ca(\cF_0)=\big\{\, (\cF_0,\cG)
\in \hil^n : \ \cG=\{g_i\}_{i\in \IN{k}} \py \|g_{i}\|^2
=a_i \ \mbox{ for }  \ i \in \IN{k}\big\}\ ,
$$ 
$$
\py 
\cS\cC_\ca(\cF_0)=\{S_\cF:\ \cF\in \cC_\ca(\cF_0)\} \inc \matpos
\ . 
$$
\een
When the initial data $(\cF_0 \coma \ca)$ are fixed, 
we shall use the notations 
$
S_0 = S_{\cF_0}$ and
$\la = \la\ua(S_0) \ . $

\pausa
We remark that we shall use the vector $\la = \la\ua(S_0)$ 
instead of $\la\da(S_0)$ for convenience (see the comments at the beginning of Section \ref{problemon}). \EOE
\end{fed}
%
%

\noindent {\bf Main problems:} (Optimal completions with prescribed norms for majorization)
 \ID \ and let $f\in \convfs$. 
\ben
\item[P1.] Give an explicit description (both spectral and geometrical) of 
$\cF\in \cC_\ca(\cF_0)$ that are the minimizers of $P_f$ 
in $\cC_\ca(\cF_0)$. 
\item[P2.] Construct a fast algorithm that efficiently computes all possible $\cF\in \cC_\ca(\cF_0)$ that are the minimizers of $P_f$ 
in $\cC_\ca(\cF_0)$. 
\item[P3.] Verify that the set of  $\cF\in \cC_\ca(\cF_0)$ that are the minimizers of $P_f$ 
in $\cC_\ca(\cF_0)$ is the same for every  $f\in \convfs$.
\een
\EOE

\pausa In previous works we have obtained some results related with the problems above. Indeed, in \cite{MRS4} we obtained a partial affirmative answer to P3, while in \cite{MRS5} we obtained some partial results related with P1. and a non-efficient algorithm as in P2. that worked in small examples (see Sections \ref{sec3} and \ref{TFC} below). 

\pausa
In this paper, building on our previous work, we completely solve the three problems above in terms of a constructive (algorithmic) approach. 


\subsection{On the structure of the minimizers of $P_f$ on $\cC_\ca(\cF_0)$}\label{sec3}
	
In this section we collect results of \cite{MRS5} that we shall use in 
this paper. Throughout this section we fix 
the initial data $(\cF_0\coma \ca)$ for the CP.
Notice that we 
are using the following convention in Definition \ref{cp}: 
we denote $\la = \la\ua(S_{\cF_0}) \in \R^d$, i.e. 
arranged in non-decreasing order. Thus we recast the 
results from \cite{MRS5} using this convention. Also notice that 
we are assuming that
$\cb=\cb^\downarrow\in \R^k$. 

\pausa
Our analysis of the completed frames $\cF= (\cF_0\coma \cG )$ depends on $\cF$ 
through $S_\cF=S_{\cF_0}+S_\cG\,$. Hence, the following description of $\cS\cC_\ca(\cF_0)$ plays a central role in our approach.

\begin{pro}\label{con el la y mayo}\rm
\ID .
Then 
\beq
\cS\cC_\ca(\cF_0) =
\big\{S\in \matpos \, :\, 
S\geq S_{\cF_0}  \py \ca \prec \lambda(S-S_{\cF_0})  
\big\}
\ .\QEDP
\eeq
\end{pro}

\pausa
Let $\mu\in (\R_{\geq 0}^d)\da$ be such that $\cb\prec \mu$, and let 
$$
\cC_\ca(\cF_0 \coma \mu)
\igdef\{ \cF=(\cF_0\coma \cG )\in \cC_\ca(\cF_0):\ \la(S_\cG)
=\mu\} 
\inc\cC_\ca(\cF_0) \ .
$$
By Proposition \ref{con el la y mayo} we get the following partition:  
\beq\label{particion}
\barr{rl}
\cC_\ca(\cF_0)& =\bigsqcup\limits_{\mu\in \,\Gamma_d(\cb)} 
\cC_\ca(\cF_0 \coma \mu)
\peso{where}  
\quad \Gamma_d(\cb)\igdef  \{ \mu\in (\R_{\geq 0}^d)\da: \ \cb\prec \mu\}
\ . \earr
\eeq
Building on Lidskii's inequality (see \cite[III.4]{Bat}) we obtained the following result:
\begin{teo}\label{dale que va} \rm 
Consider the previous notations and 
fix $\mu = \mu\da \in \Gamma_d(\cb)$. Then,
\begin{enumerate}
\item The set $\Lambda(\cC_\ca(\cF_0\coma \mu ))\igdef 
\{\lambda(S_\cF):\ \cF\in \cC_\ca(\cF_0\coma \mu )  \}$ is convex.
\item Let $\nu \igdef \la\ua(S_{\cF_0})+\mu\da $. Then $\nu\da$  is 
a $\prec$-minimizer in $\Lambda(\cC_\ca(\cF_0\coma \mu ))$.
\item If $\cF=(\cF_0\coma \cG )\in  \cC_\ca(\cF_0\coma \mu )$ 
is such that $\lambda(S_\cF)=\nu\da$ then $S_{\cF_0}$ and $S_{\cG }$ commute.
\QED
\end{enumerate}
\end{teo}

\begin{rem}\label{min OP}
Consider the previous notations and fix $\mu = \mu\da\in \Gamma_d(\cb)$. 
Let $f\in \convfs$  
and let $P_f$ be the convex potential induced by $f$.  
By the results described in Section \ref{subsec 2.2} 
and Theorem \ref{dale que va} we see that, if   $\la=\la\ua(S_{\cF_0})$ then 
\beq\label{eqpara cafop}\cF\in
{\rm argmin}\{P_f(\cF\,')  : \cF\,'\in \cC_\ca(\cF_0 \coma \mu )\} 
\iff \lambda(S_\cF)=(\lambda+\mu )\da  =(\,\la\ua +\mu \da\,)\da\ .
\eeq 
That is, if we consider the partition of $\cC_\ca(\cF_0)$ described in Eq. \eqref{particion}, 
then in each slice $\cC_\ca(\cF_0 \coma \mu )$ the minimizers of the potential $P_f$ 
are characterized by the spectral condition  \eqref{eqpara cafop}.
This shows that in order to search for global minimizers 
of $P_f$ on $\cC_\ca(\cF_0)$ we 
can restrict our attention to the set 
\beq \label{eqcafop} 
\cC_\ca^{\rm op}(\cF_0) \igdef 
\big\{\, \cF=(\cF_0\coma \cG )\in \cC_\ca(\cF_0):\ \lambda(S_{\cF})
= \big(\,\la\ua(S_{\cF_0})+\la\da(S_{\cG })\, \big)^\downarrow \, \big\} \ .
\eeq
Indeed, Eqs. \eqref{particion} and  \eqref{eqpara cafop} 
show that if $\cF$ is a minimizer of $P_f$ in $\cC_\ca(\cF_0)$ 
then $\cF\in \cC_\ca^{\rm op}(\cF_0)$.
Since the potential $P_f(\cF)$ depends on $\cF$ through the 
eigenvalues of $S_\cF$ we 
introduce the sets 
\beq\label{defi scafop} 
\cS(\cC_\ca^{\rm op}(\cF_0)) \igdef
\{S_\cF:\ \cF\in \cC_\ca^{\rm op}(\cF_0)\}
 \ \py \ \Lambda(\cC_\ca^{\rm op}(\cF_0)\,) \igdef \{\lambda(S_\cF):\ \cF\in \cC_\ca^{\rm op}(\cF_0) \}\  .
\eeq
Finally, for any $\la \in \R_{\ge0}^d \,$,  in what follows we shall also consider the set 
\beq\label{el convex}
\Lambda_\ca^{\rm op}(\lambda) \igdef
\{\lambda\ua +\mu :\mu \in \Gamma_d(\cb)\} = 
\{\lambda\ua +\mu\da : \ \mu \in \R^d_{\geq 0} \py \cb \prec \mu \}\ .
\eeq\EOE
\end{rem}

\begin{teo}\label{teocafop2}  \rm
\ID.
Denote by  $\la=\la\ua(S_{\cF_0})$. Then
\ben
\item The set $\Lambda_\ca^{\rm op}(\lambda)$ is compact and convex.
\item 
The spectral picture $\Lambda(\cC_\ca^{\rm op}(\cF_0)\,) 
=\{ \nu\,^\downarrow:\ \nu\in \Lambda_\ca^{\rm op}(\lambda)\}$.
\item If $\cF=(\cF_0 \coma \cG )\in\cC_\ca^{\rm op}(\cF_0)$, 
with $\lambda\da(S_{\cG }) =\mu$, then there exists 
$\{v_i: i\in \IN{d}\}$ an ONB
 of eigenvectors for $S_{\cF_0}\coma \la$  such that 
\beq 
S_{\cG } = \sum_{i\in \IN{d}} \, \mu_i \cdot v_i \otimes v_i 
\py 
S_{\cF}=S_{\cF_0} + S_{\cG } = \sum_{i\in \IN{d}} \, (\la_i + \mu_i) \, v_i \otimes v_i \ .
\QEDP\eeq 
\een 
\end{teo}

\pausa
For every $f\in \convf$ we consider the convex map 
\beq\label{La F}
F:\R_{\geq 0}^d \rightarrow \R \peso{given by} F(\gamma)
=\tr\,f(\ga) =
\sum_{i\in \IN{d}} f(\gamma_i) \ , \peso{for} \gamma\in \R_{\geq 0}^d\ .
\eeq
 
\begin{teo}\label{teo sobre unico espectro}  \rm
\ID \ and let $f\in \convfs$.  
Then there 
exists a vector $\elmu = \mu = \mu \da\in 
\Gamma_d(\cb)$ such that: 
\ben
\item $\cF=(\cF_0 \coma \cG )\in \cC_\ca(\cF_0)$ is a global minimizer of $P_f \iff
\cF\in \cC_\ca^{\rm op}(\cF_0)$ and $\lambda(S_{\cG })=\mu$.
\item If we let $\la = \la\ua(S_{\cF_0})$ then $\mu$ is  {\bf uniquely} 
determined by the conditions 
\beq \label{nu min LA}
\mu \in \Gamma_d(\cb) \py 
F(\la+\mu)=\min_{\gamma\in\Gamma_d(\cb)}F(\la+\gamma) =
\min_{\nu\in\Lambda_\ca^{\rm op}(\la)}F(\nu) \ . 
\eeq
Hence, if we let $\elnu \igdef \la +\elmu$ then $\exists \, ! \ \ \text{argmin}  \ \{ F(x) : {x\in \Lambda_\ca^{\rm op}(\la)}\}=\elnu$. 
\qed
\een
\end{teo}


\begin{teo}\label{vale la conjetur} \rm
\ID \,. Let $f\in \convfs$  
and  assume that $\cF=(\cF_0\coma \cG)$ is a global 
minimizer of $P_f$ on $\cC_\ca^{\rm op}(\cF_0)$. 
Then, there exists a partition $\{J_i\}_{i\in \IN{p}}$ of $\IN{k}\,$ and $c_1>\ldots>c_p>0$ 
such that 
\ben
\item The subfamilies $\cG_i=\{f_j\}_{j\in J_i}$ (for $i\in\IN{p}\,$) are mutually orthogonal, 
i.e. $S_\cG = \oplus_{i\in\IN{p}} S_{\cG_i}\,$.
\item
The frame operators  $S_{\cG_i}$ and $S_{\cF_0}$ commute, for every $i\in \IN{p}$ . 
\item
We have that 
$ 
S_\cF \, f_j=c_i\, f_j\, $, for every 
$j\in J_i$ and every $i\in \IN{p}$ .
\end{enumerate}
The statement is still valid if we assume that $\cF$ is just a local minimizer, 
but if we also assume as a hypothesis that $\cF$  satisfies item 2 
(for example if $S_{\cF_0}=0$). \QED
\end{teo}

\subsection{The feasible case of the CP}\label{TFC}
\pausa
In this section we recall the results from \cite{MRS4} that we shall need in the sequel. Throughout this section we fix
the initial data $(\cF_0,\ca)$ for the CP.
Denote by $S_0 = S_{\cF_0}\,$, $\la = \la\ua (S_0)$ and 
$t=\tr\, \la + \tr \, \ca$. In \cite{MRS4} we introduced the following set 
$$ 
U_t(S_0\coma m)=\{S_0+B:\   B\in \matpos \, , \ \rk \, B 
\le d-m \ , \ \tr\,(S_0+B)\ =\  t\ \}\inc\matpos \ ,
$$
where $m = d-k$. In \cite[Theorem 3.12]{MRS4} 
 it is shown that there exist $\prec$-minimizers in $U_t(S_{0}\coma m)$.  
Indeed, there exists 
$\mu (\la \coma \ca )\in (\R^d_{\geq 0})\da$ - that 
can be effectively computed by a fast algorithm - such that, 
if $\nu (\la \coma \ca )\igdef \la+\muel\in \R^d_{>0}$ 
then $S\in U_t(S_{0}\coma m)$ is  a $\prec$-minimizer 
if and only if $\la(S)=\nu (\la \coma \ca )\da$.

\pausa
Notice that by construction $\nu (\la \coma \ca )$ is not a necessarily ordered vector 
(nor decreasing, nor increasing); yet, in terms of the terminology from 
\cite{MRS4}, we have that $\nu_{\la\coma m}(t)=\nu (\la \coma \ca )\da$. Thus, we have reversed the order of the vector $\muel$ 
- accordingly with reversing the order of $\la=\la\ua(S_{\cF_0})$ -  and we have changed the description of the vector $\nuel$ - while preserving all of their majorization properties - with respect to \cite{MRS4}. Nevertheless, we point out that the ordering of the entries of the vector $\nu (\la \coma \ca )$ presented here plays a crucial role in simplifying the exposition of the results herein, as it guaranties that $\muel =\nuel-\la$. 

\pausa
The following definition and remark show the relevance of the notions introduced above for the computation of the spectral structure of solutions for the optimal 
 completion problem.
\begin{fed}\label{defi fea}\rm 
\ID \ with $\la = \la\ua(S_{\cF_0})$. We say that the pair 
$(\la\coma \ca)$ is {\bf feasible} if 
$\muel$ satisfies that $\cb\prec\muel$. \EOE
\end{fed} 

\begin{rem}\label{rem fea}
\ID \ with $\la = \la\ua(S_{\cF_0})$. Assume that the pair $(\la\coma \ca)$ is feasible and denote $\mu=\muel$. 
In this case (see \cite{MRS4}) for any $S$ which is a $\prec$-minimizer 
in $U_t(S_0\coma m)$ - where $m=d-k$ - it holds that $\la(S-S_0)=\mu$ and 
hence, by Proposition \ref{con el la y mayo}, we conclude 
that $S\in \cS\cC_\ca(\cF_0)$. Moreover, Proposition 
\ref{con el la y mayo} also shows that  $\cS\cC_\ca(\cF_0)\inc U_t(S_0\coma m)$.  
Then $S$ is also a $\prec$-minimizer in $\cS\cC_\ca(\cF_0)$. 
Therefore, as a consequence of the results in 
Section \ref{subsec 2.2}, 
any completion $\cF=(\cF_0,\cG)\in \cC_\ca(\cF_0)$ such that $S_\cF=S$ is a 
minimizer of $P_f$ for every $f\in \convf$. 
  
\pausa
On the other hand, as a consequence of the geometrical structure of $S=S_\cF$ as above 
(see \cite{MRS4,MRS5}), we conclude that there exists $c>0$ such that $S_\cF \,g_i=c\, g_i$ 
for every $i\in \IN{k}\,$. That is, in this case the structure of the completing sequence 
$\cG$ given in Theorem \ref{vale la conjetur} is trivial: the partition of $\IN{k}$ has only one member and there exists a unique constant $c=c_1\,$. 
\EOE
\end{rem}
\pausa
It is worth pointing out that it is easy to construct examples of initial data $(\cF_0\coma
\ca)$ for the CP such that the pair $(\la,\ca)$ is not feasible (see \cite{MRS4}), so that comments in Remark \ref{rem fea} do not apply in these cases.

\begin{rem}\label{nuel con m0}
\ID \ with $k\geq d$ (so that $m=d-k\leq 0$), let 
$\la = \la\ua(S_{\cF_0})$ and let $t=\tr\, \ca+\tr \, \la$. 
In \cite{MRS4} we shown that there are two cases: 
\ben
\item Since $\la = \la\ua$ then $\la_d 
= \max\,\{ \la_i : i \in \IN{d}\}\,$. If
\beq\label{nu la ca truc}
\frac td = \frac{\tr\, \ca + \tr\,\la }{d} 
\geq \la_d  \peso{then}  \la \leqp \frac{t}{d}  
\ \uno_d = \nuel \ .
\eeq
\item If $\la_d > \frac td$ then there exists $s\in \IN{d-1}$ such that 
\beq\label{nu la ca}
\nu (\la \coma \ca ) 
= (c\, \uno _{s}\coma  \la_{s+1} \coma \dots \coma \la_d ) 
\peso{with} \la_{s}\le c< \la_{s+1} 
\eeq
so that $\la \leqp \nuel = \nuel\ua$,  and in this case the index $s$
also satisfies that 
$$
c = \frac{1}{s} \ \big[\  \tr\, \ca + \sum_{i=1}^{s}\  \la_i 
\, \big] \peso{so that} \tr\, \nuel = t 
= \tr\, \la + \tr \, \ca  \ . 
$$ 
\een
In what follows we obtain an explicit description of the 
vector $\nuel$ in case $d\leq k$ (so that $m\leq 0$) and 
$\frac{1}{d}\, [\, \tr\, \ca + \tr\,\la \, ]< \la_d\,$. 
Explicitly, we compute the parameters $s$ and $c$ of 
Eq. \eqref{nu la ca} in a way that is key for the developments 
of Section \ref{sec The main result}. Our present techniques 
differ substantially from those introduced in \cite{MRS4}.
%
We begin by showing that the
vector $\nuel$ above is unique. Then, we show that the computation 
of $\nuel$ for $m\in\IN{d-1}$ can be reduced to the case when $m=0$.
First we need to introduce some notations: 
\EOE
\end{rem}

\begin{fed}\rm
\ID . Assume that $d\le k$. 
We denote by  
$$
\la = \la\ua(S_{\cF_0}) \py 
h_i=\la_{i}+a_i \peso{for every} i \in \IN{d} \ .
$$ 
Given  $j\coma r \in \IN{d}\cup\trivial$ such that $ j< r$, 
by $Q_{j\coma r}$ we denote the final averages:
\beq\label{los Q bis}
Q_{j\coma r} =\frac{1}{r-j} \ \Big[\  \sum_{i=j+1}^{r}\  h_i 
+ \sum_{i=r+1}^k\  a_i  \Big] 
=\frac{1}{r-j} \ \Big[\ \sum_{i=j+1}^{k}\  a_i + 
\sum_{i=j+1}^{r}\  \la_i    \Big] \ .
\eeq
We shall abbreviate $Q_r = Q_{0\coma r}\,$. \EOE
\end{fed}

\begin{lem} \label{el r de M bis}
\ID \ with $k\ge d$. Let $r \in \IN{d}\,$.  Then 
\ben
\item If $r<d$ and $Q_r<\la_{r+1}$ then $Q_r<Q_j\,$, for every $j$ such that $r<j\le d$. 
\item If $r<d$ and $Q_r\le\la_{r+1}$ then $Q_r\le Q_j\,$, for every $j$ such that $r<j\le d$. 
\item If $\la_r\le Q_r$ then $Q_r\le Q_j\,$, for every $j$ such that $1\le j< r$. 
\een
\end{lem}
\proof 
Denote by $c= Q_r\,$ for a fixed $r<d$. Recall that $\la = \la\ua$. If  $j>r$ then 
\begin{align*}
c < \la_{r+1} \implies 
Q_j& =\frac{1}{j}\left(\,\tr\, \ca + \sum_{i=1}^r \la_i +
\sum_{i=r+1}^j \la_{i}\right) 
> \frac{1}{j}\ (\,r\,c+(j-r)\,c\,)=c \ .
\end{align*}
The proof of item 2 is identical. On the other side, if $j< r$ then 
\begin{align*}
\la_{r}\leq c \implies 
Q_j& 
=\frac{1}{j}\left(\,\tr\, \ca + \sum_{i=1}^r \la_i -
\sum_{i=j+1}^r \la_{i}\right) 
\geq \frac{1}{j}\ (\, r\, c-(r-j)\, c\, )=c \ . \QEDP
\end{align*}

\begin{pro} \label{el r de M tris}
\ID \ with $k\ge d$ (so that $m\le 0$) and assume that 
$\frac{1}{d}\, [\, \tr\, \ca + \tr\,\la 
\, ]< \la_d$. Then 
\ben
\item There exists a unique  index
$s\in \IN{d} $ such that $\la_{s}\le  Q_{s}<\la_{s+1}\,$, and in this case
\beq\label{EL r}
s = \max \ \{ w\in \IN{d-1}:  Q_{w} = 
\min\limits_{j\in \IN{d}}\  Q_{j} \ \} \peso{and} \nu(\la \coma  \ca)  
= (Q_s\, \uno _{s}\coma  \la_{s+1} \coma \dots \coma \la_d )
\eeq
\item If another index $r \in \IN{d-1}$ satisfies that 
$\la_{r}\le  Q_{r}\le \la_{r+1}\,$, then 
\ben
\item  $Q_{r} = \min\limits_{j\in \IN{d}}\  Q_{j} = Q_s\,$ and $ r\le s$. 
\item  If  $r<s$, then  $Q_r = \la_{r+1} = 
\la_s\,$ and also $\nuel=(Q_r\, \uno _{r}\coma  \la_{r+1} \coma 
\dots \coma \la_d ) $.
\een
\item Given $\rho = (c\, \uno _{r}\coma  \la_{r+1} \coma 
\dots \coma \la_d )$ (or $\rho = c\, \uno_d$) such that $\la \leqp \rho = \rho\ua$ and $\tr \, \rho = 
\tr \, \nuel$ then 
$\rho = \nuel$. 
\een
\end{pro}
\proof The existence of an index $s$ such as in item 1 is guaranteed by 
the properties of $\nuel$ stated in \cite{MRS4}. Nevertheless, 
it is easy to see that the index $s$ described in Eq. \eqref{EL r} 
satisfies that $\la_{s}\le  Q_{s}<\la_{s+1}\,$. The formula given in Eq. \eqref{EL r}, which shows 
the uniqueness of $\nuel$, is a direct consequence of Lemma \ref{el r de M bis}.
Assume that $\la_{r}\le  Q_{r}\le \la_{r+1}\,$. Then $Q_{r} = \min\limits_{j\in \IN{d}}\  Q_{j} = Q_s\,$ 
and $ r\le s$ by Lemma \ref{el r de M bis}.  
If $r<s$, then $Q_{s} = \frac1{s} \, (r\, Q_r + \sum_{i=r+1}^s\la_i ) = Q_r\,$. 
This clearly implies all the equalities of item (b). 
Finally, observe that item 2 $\implies$ item 3. \QED

\begin{rem}\label{cuando m>0}
\ID \ with $m = d-k >0$. Then if
\beq\label{el nu la a}
\tilde\la  = (\la_{1}\coma \dots \coma \la_{k}) \in (\R^{k})\ua
\peso{then}
\nu (\la \coma \ca) = (\nu (\,\tilde\la \coma \ca) \coma  
\la_{k+1} \coma  \ldots \coma  \la_d)\ ,
\eeq
and $\nu (\, \tilde\la \coma \ca)$ is constructed as in 
Proposition \ref{el r de M tris}.  

\pausa
The proof is direct by 
observing that, 
extracting the entries $\la_{k+1} \coma  \ldots \coma  \la_d$ 
of the vector $\nuel$ as described in \cite[Def. 4.13]{MRS4},  
the vector that one obtains  
(with the reverse order) 
satisfies the conditions 
of item 3 of Proposition \ref{el r de M tris} relative to the pair 
$(\,\tilde\la \coma \ca)$.  
\EOE	
\end{rem}

\pausa
The following result is in a sense a converse to Remark \ref{rem fea}. It establishes that if there exists $f\in \convfs$ 
and a minimizer $\cF=(\cF_0,\cG)$ of $P_f$ in $\cC_\ca(\cF_0)$ such that the structure of the completing sequence $\cG$ as described in Theorem \ref{vale la conjetur} is trivial, the the underlying pair $(\la\coma \ca)$ is feasible. Recall the notation 
$\elnu $ given in Theorem \ref{teo sobre unico espectro}.

\begin{lem}\label{solo c1} \ID \,, $k\ge d $, 
and let $\cF= (\cF_0\coma \cG)\in \cC_\ca^{\rm op}(\cF_0)$ 
be a  minimum for $P_f\,$ on $\cC_\ca(\cF_0)$ for a $f\in \convfs$. 
Suppose that, for some  $c >0$, 
$$
 W=R(S_\cG)\neq  \cH   \py  S_\cF\big|_W \in L(W)= c\,I_W\ .    
$$
Let $\la = \la\ua(S_{\cF_0})$, $\mu = \la\da(S_\cG)$ and 
$s\igdef \dim \, W = 
\max \{i \in \IN{d}: \mu_i\neq0\}$. Then 
$$
\la_s < c \le \la_{s+1} \peso{so that}
(\la\coma \ca) \peso{is feasible and}
\elnu  =\nuel \ .
$$
The same final conclusion trivially
holds if $s=\dim \, W = d$ and $S_\cF = c \, I$. 
\end{lem}
\proof 
Suppose that $s<d$. By hypothesis  $\elnu 
=\la\ua+\mu\da=\big(\,c\, \uno_s
\coma   \la_{s+1}  \coma \dots \coma \la_d \big) $ 
and it satisfies that $\la(S_\cF)=\elnu\da$.  
Since 
$\ca\prec \mu = \mu\da$ then $
 \tr \, \mu = \tr \, \ca > \sum_{i=1}^{s} \, a_i \,$, 
because $s < d\le k $. 
Suppose now that $c  >\la_{s+1}\,$. For small $t>0$
consider the vector 
$$
\ga(t) = \big(\,c \, \uno_{s-1} \coma (c-t) \coma   
\la_{s+1}+t \coma \la_{s+2}\coma 
\dots  \coma \la_d   \big) \in \R^d \peso{with}  
\tr \,\ga(t) = \tr\,S_\cF\ . 
$$ 
Let $\mu(t) = \ga(t)-\la$. 
For every $t$ we have that   $\tr \, \mu(t) = \tr\,\mu$. 
On the other hand, if
$$
t<\frac{\mu_s}2 \implies 
\mu(t) = (\mu_1\coma \dots \coma \mu_{s-1}
\coma \mu_s-t \coma t  \coma 0\, \uno_{d-s-1}) = \mu(t)\da 
\in (\R_{\ge0}^d)\da\ .
$$
It is easy to see that  
if also $t <\sum_{i=s+1}^k \, a_i $ then still $\ca \prec \mu(t)$. 
So there exists 
$\cF'
\in \cC_\ca^{\rm op}(\cF_0)$ such that
$\la(S_{\cF'}) = \ga(t)\da$. 
Notice that, since  $(c-t  \coma \la_{s+1}+t) \prec 
(c\coma \la_{s+1})$ strictly, then 
$P_f(\cF') = \tr \, f(\ga(t)\,) <\tr\, f(\elnu\,) = P_f(\cF)$, 
a contradiction. Hence $c\le \la_{s+1}\,$. 

\pausa
The condition $\la_s< c$ follows from the fact that 
$ c-\la_s=\mu_s>0$. 
These facts show that  $\la= \la\ua\leqp \elnu = \elnu\ua \implies 
\elnu= \nu(\la\coma \ca)$ (by item 3 
of Proposition \ref{el r de M tris}). 
In particular, $\ca\prec \la(S_\cG) = \mu = \nuel-\la =\muel$ 
so that $(\la\coma \ca)$ is feasible. \QED

\section{Uniqueness and characterization of the minimum}\label{sec The main result}
In this section we shall state the main results of th paper. For 
the sake of clarity of the exposition,  we postpone the more technical proofs
until Section 5. 

\begin{num}[Fixed data, notations and terminolgy]\label{notaciones conj}
\ID. Until  Theorem \ref{final remark}, we 
shall assume that $k\ge d$, so that $m=d-k \le 0$. Recall that 
$$\cC_\ca^{\rm op}(\cF_0) \igdef 
\big\{\, \cF=(\cF_0\coma \cG )\in \cC_\ca(\cF_0):\ \lambda(S_{\cF})
= \big(\,\la\ua(S_{\cF_0})+\la\da(S_{\cG })\, \big)\da \, \big\} \ .
$$ 
Fix $f\in \convfs$ and 
a minimizer $\cF= (\cF_0\coma \cG)\in \cC_\ca(\cF_0)$ be 
 for $P_f\,$ on $\cC_\ca(\cF_0)$.
\ben
\item  By Theorem \ref{teo sobre unico espectro}, we know that 
$\cF\in \cC_\ca^{\rm op}(\cF_0)$ and, if we denote by 
$\la = \la\ua(S_{\cF_0})$, then 
$\la\da(S_{\cG })=\elmu = \elnu-\la$.
By  Theorem \ref{teocafop2}
there exists $\{v_i: i\in \IN{d}\}$ an ONB
of eigenvectors for $S_{\cF_0}\coma \la$  such that 
\beq\label{la BON2bis}
S_{\cG } = \sum_{i\in \IN{d}} \, \mu_i \cdot v_i \otimes v_i 
\py 
S_{\cF}=S_{\cF_0} + S_{\cG } = \sum_{i\in \IN{d}} \, (\la_i + \mu_i) \, v_i \otimes v_i \ .
\eeq
\item Let $ s_\cF= \max \, \{i \in \IN{d}
: \mu_i \neq 0\} = \rk \, S_\cG\,$. Denote by  $W= R(S_\cG)$, which reduces $S_\cF\,$.  
\item Let $S = S_\cF\big|_W \in L(W)\,$ and $\sigma(S) = \{ c_1 \coma \dots \coma c_p\}$
(where $c_1 > c_2 > \dots > c_p>0$). 
\item Let $K_j = \{ i \in   \IN{s_\cF}:  \la_i +\mu_i = c_j\}$  
and $J_j = \{i \in \IN{k}: S\,g_i = c_j \, g_i\}$. 
By Theorem \ref{vale la conjetur},  
$$
\barr{rl}
\IN{s_\cF} & = \bigsqcup\limits_{j\in \IN{p}} 
\, K_j \py 
\IN{k} = \bigsqcup\limits_{j\in \IN{p}}
\, J_j  \ . \earr
$$
Observe that $s_\cF\,$, $\la = \la\ua(S_{\cF_0})$ and the sets $K_j$ completely describe the vector
$\mu = \la(S_\cG)$.
\item Since $R(S_\cG)= \gen\{g_i : i \in \IN{k}\} = W = \oplus_{i \in \IN{p}} \ker \,(S-c_i\,I_W\,)$
then for every $j\in \IN{p}\,$, 
\beq\label{cajas}
W_j \igdef \gen\{g_i : i \in J_j\} = \ker \,(S-c_j\,I_W\,) = \gen\{v_i : i \in K_j\} \ ,
\eeq
because $g_i \in  \ker \,(S-c_j\,I_W\,)$ for every $i \in J_j\,$. Note that, by Theorem \ref{vale la conjetur}, each $W_j$ reduces both $S_{\cF_0}$ and $S_\cG\,$. 
\item If $p = 1$ then $J_1 = \IN{k}\,$ and $S = c_1\, I_W\,$.  
Hence the minimum $\cF$ satisfies the hypothesis of Lemma \ref{solo c1}, 
so that the pair $(\la\coma \ca) $ is feasible. 

\item We denote by  $h_i=\la_{i}+a_i$ for every $i \in \IN{d}\,$. 
Given  $j\coma r \in \IN{d}$ such that $ j\le r$, let 
\[
P_{j\coma r} =\frac{1}{r-j+1}\ \sum_{i=j}^r\  h_i = 
\frac{1}{r-j+1}\ \sum_{i=j}^r\ \la_{i}+a_i 
\ ,
\] 
be the initial averages. We abbreviate $P_{1 \coma r} = P_r\,$. \EOE  
\een
\end{num}

\begin{rem}[A reduction procedure] \label{induc}\rm
Consider the data, notations and terminology fixed 
in \ref{notaciones conj}. 
For any $j \in \IN{p-1}$ denote by 
$$
I_j = \IN{d} \setminus \bigcup_{i \le j} \,K_i  \ \coma  \ L_j = 
\IN{k} \setminus \bigcup_{i \le j} \,J_i  \ \coma \ 
\la^{(j)} = (\la_i)_{i \in I_j}  \ \coma  \ 
\cG_j = (g_i)_{i \in L_j}  \coma \ca^{(j)}= (a_i )_{i \in L_j} 
$$
and take some sequence $\cF_0^{(j)} $ in 
$\cH_j = \big[\, \bigoplus_{i\le j} W_i \big]\orto$
 such that $S_{\cF_0^{(j)} } = S_0|_{\cH_j}$ (notice that, by construction, $\cH_j$ reduces $S_0$). 

\pausa
Then, it is straightforward to show that $\cF_j = (\cF_0^{(j)} \coma \cG_j) $ is a (global) minimizer of 
$P_f$ on $\cC_{\ca_j}(\cF_0^{(j)})$ in $\cH_j\,$, i.e. an optimal completion for the reduced problem. 
Indeed, recall that the minimality 
is computed in terms of the map $F$ defined in Eq. \eqref{La F}, 
which works independently in each entry 
of $\lambda(S_\cF)=\elnu\da$. 

\pausa
The importance of the previous remarks lies in the fact that 
they provide a powerful reduction method to compute the structure 
of the sets $\cG_i \coma K_i $ and $J_i$ for $i\in \IN{p}$ as well as 
the set of constants $c_1>\ldots>c_p>0$. Indeed, assume that we are 
able to describe the sets $\cG_1\coma K_1\coma J_1$ and the constant 
$c_1$ in some structural sense, using the fact that these sets are 
extremal (e.g. these sets are built on $c_1>c_j$ for $2\leq j\leq p$).

\pausa
Then, in principle, we could apply these structural arguments to 
find $\cG_2\coma K_2\coma J_2$ and the constant $c_2\,$, using the fact that 
these are now extremal sets of $\cF_1\,$, which is a $P_f$ minimizer 
of the reduced CP for $(\cF_0^{(1)},\ca^{(1)})$. 
On the other hand, the minimality of the final reduction  $\cF_{p-1}\,$ 
produces a pair $(\la^{(p-1)}\coma \ca^{(p-1)})$ which must 
be feasible by item 6 of \ref{notaciones conj}, 
because it has an unique constant $c_p\,$ associated to the unique set $K_p\,$. 
As we shall see, 
this strategy can be implemented to obtain (inductively) a 
precise description of the sets above.
 \EOE
\end{rem}

\begin{rem} \label{remando0}
\ID \ with $d\le k$, $\la = \la\ua(S_{\cF_0})$ and $\cb= \ca\da$.
Fix $f\in \convfs$ and let 
$\cF= (\cF_0\coma \cG)\in \cC_\ca^{\rm op}(\cF_0)$ be a { global} minimum for $P_f\,$
on $\cC_\ca^{\rm op}(\cF_0)$.  
In section \ref{proof remando} we shall prove the following properties (conjectured in \cite{MRS5})
of the sets $J_j$ and $K_j$ defined in item 4. of \ref{notaciones conj} describing $\elmu$ and $\elnu$: 
%
\ben 
\item Each set $J_j$ and $K_j$  
consists of { consecutive} indexes, for $j\in\IN{p}$ .  
\item The sets $K_j$ and $J_j$ have the same number of elements, for $j\in \IN{p-1}\,$.

\item Moreover, $J_1<\ldots<J_p$ (i.e. if $l\in J_i$ and $h\in J_j$ with $i<j \, \Rightarrow \, l<h$) and $K_1<\ldots<K_p$. In particular, by items 1 and 2 above, $K_j=J_j$ for $j\in \IN{p-1}\,$.
 \een
\EOE
\end{rem}

\pausa
We state the properties of the sets $J_j$ and $K_j$, $j\in\IN{p}$ described in Remark \ref{remando0} in the following:

\begin{teo}\label{remando} \rm
\ID \ with $d\le k$. With the notations of Remark \ref{remando0}, 
assume that $\la = \la\ua(S_{\cF_0})$, $\mu = \mu\da = \elmu$ and $\ca = \ca\da$.  
Then 
\ben
\item 
There exist 
$0 = s_0<s_1<  s_2<\cdots< s_{p-1}<s_p = s_\cF= \max \{j \in \IN{d} : \mu_j \neq 0\}$
such that 
\[
K_j =J_j =\{s_{j-1}+1\coma \ldots \coma  s_j\} \ ,\quad 
\ ,\peso{for} j\in\IN{p-1}\  , 
\]
\[
 K_p=\{s_{p-1}+1 \coma \ldots \coma  s_p\} \ , \ J_p=\{s_{p-1}+1\coma \ldots\coma  k\} \ .
\] 
\item The vector $\nu_f(\la\coma \ca) = \big(\, c_1 \, \uno_{s_1} \coma 
\dots \coma c_p\, \uno_{s_p-s_{p-1}} \coma \la_{s_p+1} \coma 
\dots \coma \la_d\,\big)$, where 
$$
c_r=\frac{1}{s_r-s_{r-1}} \ \sum_{i=s_{r-1}+1}^{s_r}\, h_i
= P_{s_{r-1}+1\coma s_{r}}\peso{for} r\in\IN{p-1}\ ,
$$
or also $c_r = \la_j + \mu_j$ for every $j \in K_r=J_r $ for $ r\in\IN{p-1}\, $. 
\item The constant $c_p$ is the one defined by the feasible final part i.e.,  
 $c_p = Q_{s_{p-1}\coma  s_p}$ and the indexes $s_{p-1}$ and $s_p$ are determined by 
the last block  (recall  Lemma \ref{solo c1}). 
\een
\end{teo}
\proof See Section \ref{sect sev proof}.
 \QED

\pausa
\ID . Assume that $\elnu=\big(\, c_1 \, \uno_{s_1} \coma 
\dots \coma \la_{s_1+1} \coma 
\dots \coma \la_d\,\big)$ i.e. with $p=1$, in the notations of 
Theorem \ref{remando}. Then, by Lemma \ref{solo c1}, 
the pair $(\la \coma \ca)$ is feasible and $\elnu = \nuel$. 

\pausa
In what follows we shall need the following notion, that allow us to show feasibility in the more general case in which, in the notations of 
Theorem \ref{remando}, $p>1$.

%
%
%

\begin{fed} \label{s feas}\rm 
\ID . Let $\la = \la\ua(S_{\cF_0})\in (\R_{>0}^d)\ua$. Suppose that $k\ge d$. 
Given $s\in \IN{d-1}$ denote by 
$$\la ^s = (\la_{s+1} \coma \dots \coma 
\la_d )\in \R^{d-s}  \py \ca^s = (a_{s+1}\coma \dots \coma a_k) 
\in \R^{k-s}\ ,
$$
the truncations of the original vectors $\la$ and $\ca$. 
We say that the index $s$ is {\bf feasible} if the pair 
$(\la^s\coma \ca^s)$ is feasible for the CP. 
Note that $(d-s)-(k-s) = d-k= m\le 0$. Therefore 
$$
\nu_s \igdef  \nu (\la^s\coma \ca^s) \stackrel{\eqref{nu la ca}}{=} 
\big(\,c\, \uno_{r-s}\coma   \la_{r+1}  \coma \dots \coma \la_d  \big)
\peso{where} c = Q_{s \coma r} 
$$
for the unique $r> s$ such that $\la_{r}\le c<\la_{r+1}\,$ (or 
$\nu_s = Q_{s\coma d}\, \uno_{d-s}$ 
if $\la_d\le Q_{s\coma d}$). 
This means that $\la_s\leqp \nu_s\in (\R^{d-s}_{_{>0}})\ua$. 
\EOE 
\end{fed}

\pausa
By Remark \ref{induc} and Lemma \ref{solo c1}  we know 
that, with the notations of \ref{notaciones conj},  the index $s_{p-1}$ 
associated to the minimum $\nu= \elnu$ is feasible - in the sense of Definition \ref{s feas} -  because 
the last block of $\nu$ is constructed with the final feasible 
parts of $\la$ and $\ca$, and 
$\nu_{s_{p-1}}=
\big(\,c_p\, \uno_{s_p-s_{p-1}}\coma   \la_{s_{p}+1}  \coma \dots \coma \la_d  \big) $. 

\begin{pro}\label{s = sf bis}
\ID . 
With the notations of Theorem \ref{remando}, the global minimum $\elnu$ satisfies that 
\ben
\item The index $s_{p-1}$ (where the feasible part begins) is determined by 
$$
 s_{p-1} 
= \min \ \{\ s\in \IN{d} \ : \  s  \ \mbox{ \rm is feasible } \} \ .
$$
\item 
The following recursive method 
allow to describe the vector $\elnu$ as in Theorem \ref{remando}: 
\ben
\item The index $s_1 = \max \, \big\{j \le s_{p-1} \, :\, 
P_{1\coma j} = \max\limits_{i\le s_{p-1}}  \, P_{1\coma i} \, \big\}$, and 
$c_1 = P_{1\coma s_1}\,$.
\item If the index $s_j$ is already computed and $s_j<s_{p-1}\,$, then
\een
\een
$$
s_{j+1} = \max \, \big\{s_j< r \le s_{p-1} \, :\, 
P_{s_j+1\coma j} = \max\limits_{s_j< i\le s_{p-1}}  \, P_{s_j+1 \coma i} \, \big\} 
\py c_{j+1} = P_{s_j+1\coma s_{j+1}}\ .
$$
\end{pro}
\proof 
See Propositions \ref{s = sf} and \ref{el s_1}. \QED

\pausa
The following are the main results of the paper. In order to state them, we introduce the spectral picture of the completions with prescribed norms, given by
$$ \Lambda(\cC_\ca(\cF_0)\,) \igdef\{\la(S_\cF):\ \cF\in\cC_\cb(\cF_0)\}\ .$$

\begin{teo} \label{el teo}
\ID \ with $m= d-k\le 0$. Then the vector  $\nu = \elnu$ is the same for every 
$f\in \convfs$. Therefore,  
\beq\label{con mayo}
\nu\da \in \Lambda(\cC_\ca(\cF_0)\,) \py 
\nu\da \prec \ga \peso{for every} \ga\in \Lambda(\cC_\ca(\cF_0)\,) \ .
\eeq
\end{teo}
\proof
By Proposition  \ref{s = sf bis}, the minima 
$\nu = \elnu$ are completely characterized by the data $(\la\coma \ca)$ 
without interference of the map $f$. Therefore, 
given any $\ga\in \Lambda(\cC_\ca(\cF_0))$, 
\beq
\tr \, f(\nu) \le \tr \, f(\ga) \peso{for every}
f\in \convfs  \implies \nu \prec \ga\ . \QEDP
\eeq

\pausa
The following result shows that the structure of optimal completions in $\cC_\ca(\cF_0)$ in case $m=d-k>0$ can be obtained from the case in which $m=0$.
\begin{teo}\label{final remark}
\ID \ with $ m = d-k>0$.  If we let
$$
\la' = (\la_1\coma \dots \coma \la_{k}) \in (\R_{\geq 0}^{k})\ua 
\peso{then} \nu_f (\la \coma \ca) 
= (\nu_f (\la' \coma \ca) \coma  \la_{k+1} \coma  \ldots \coma  \la_d)\ ,
$$ 
where $\nu_f (\la' \coma \ca)$ is constructed as in 
Proposition \ref{s = sf bis} (since $d\,'= k$, by construction of $\la '\in(\R_{\geq 0}^{d\,'})\ua$). 
In this case 
the vector $\elnu$ is the same for every 
$f\in \convfs$ and
also satisfies Eq. \eqref{con mayo}. 
\end{teo}
\proof
Since $k = d-m$ and  $\ca\in \R^k$  
we deduce that any $\delta = \delta\da\in \R^d_{\ge 0}$ 
such that $\ca\prec \delta$ 
must have $\delta_{k+1} =  \ldots = \delta_d = 0$. 
It is easy to see that this fact implies that 
\beq\label{con ceros}
\Lambda_\ca^{\rm op}(\lambda) = 
\{\lambda\ua +\delta\da : \ \delta \in \R^d_{\geq 0} 
\py \cb \prec \delta \} = 
\{(\ga \coma  \la_{k+1} \coma  \ldots \coma  \la_d) : \ga \in \Lambda_\ca^{\rm op}(\la')\}
\eeq
We know that $\elnu -\la = \mu =  \mu \da$ and that 
$\ca\prec \mu \implies  
\mu_{k+1} =  \ldots = \mu_d = 0$. Recall the map $F:\R_{\geq 0}^d \rightarrow \R $ 
defined in Eq. \eqref{La F} for each $f \in \convfs$. Therefore 
\beq\label{ga ro}
\elnu\in \Lambda_\ca^{\rm op}(\lambda) \py
\nu_f (\la \coma \ca) = \mu\da+\la\ua \implies 
\nu_f (\la \coma \ca) = (\rho \coma  \la_{k+1} \coma  \ldots \coma  \la_d)\ ,
\eeq
for some $\rho \in \Lambda_\ca^{\rm op}(\la')$. Then 
$F(\elnu\,) = F(\rho) + F(\la_{k+1} \coma  \ldots \coma  \la_d)$. 
By  Eq. \eqref{nu min LA}, 
$$
F(\elnu\,)  
\stackrel{\eqref{nu min LA}}{=}  
\min_{\nu\in\Lambda_\ca^{\rm op}(\la)}F(\nu)  
\stackrel{\eqref{con ceros}}{=}  
\left[\, \min_{\ga\in\Lambda_\ca^{\rm op}(\la')} \, F(\ga) \, 
\right]+ F(\la_{k+1} \coma  \ldots \coma  \la_d)  \ .
$$
Using  Eq. \eqref{nu min LA} again we deduce that 
$\rho = \nu_f (\la' \coma \ca)$. 
Since $\nu_f (\la' \coma \ca)$ 
is constructed as in 
Proposition \ref{s = sf bis}, then it 
is the same vector for every strictly convex map $f$ and  the same 
happens with $\nu_f (\la \coma \ca)$, so that 
$\elnu\da$  is a minimum for majorization on $\Lambda(\cC_\ca(\cF_0)\,)$. 
\QED

\begin{rem} The construction of the minimum $\elnu$ 
given by Proposition \ref{s = sf bis} is 
algorithmic, an it can be easily implemented in Mathlab.  It only depends on - an already available, see \cite{MRS4} - routine for checking feasibility, which is fast and efficient. 
\EOE
\end{rem}

\section{Proofs of some technical results.}
In this section we present detailed proofs of several statements in section \ref{sec The main result}. All these 
results assume that the initial data $(\cF_0\coma \ca)$ for the CP satisfies 
that $k\ge d$. As already explained, the general case can be reduced to this situation.

\subsection{Description of the sets $K_i$ and $J_i\,$.}\label{proof remando} 
\begin{num}
\label{notaciones conj loc}
We begin by recalling the notations of \ref{notaciones conj}:
\ID , with $k\ge d$. Fix a convex map $f \in \convfs$. We consider the following objects: 
\ben
\item Let $\cF= (\cF_0\coma \cG)\in \cC_\ca^{\rm op}(\cF_0)$ be a 
global minimum for $P_f\,$
on $\cC_\ca^{\rm op}(\cF_0)$ (or a local minimum if $\cF_0 = \vacio $).  
If $\la = \la\ua(S_{\cF_0}) $ and $\mu\igdef \lambda\da(S_{\cG })$, then there exists 
$\{v_i: i\in \IN{d}\}$ an ONB
of eigenvectors for $S_{\cF_0}\coma \la$  such that 
$$
S_{\cG } = \sum_{i\in \IN{d}} \, \mu_i \cdot v_i \otimes v_i 
\py 
S_{\cF}=S_{\cF_0} + S_{\cG } = \sum_{i\in \IN{d}} \, (\la_i + \mu_i) \, v_i \otimes v_i \ .
$$
\item Let $s_\cF= \max \{i \in \IN{d} : \mu_i \neq 0\} = \rk\,S_\cG\,$. 
Denote by  $W= R(S_\cG)$, which reduces $S_\cF\,$.  
\item Let $S = S_\cF\big|_W \in L(W)\,$ and $\sigma(S) = \{ c_1 \coma \dots \coma c_p\}$
(where $c_1 > c_2 > \dots > c_p$). 
\item Let $K_j = \{ i \in   \IN{s}:  \la_i +\mu_i = c_j\}$  
and $J_j = \{i \in \IN{k}: S\,g_i = c_j \, g_i\}$. 
Then 
$$\barr{rl}
\IN{s_\cF} &= \stackrel{_{_D}}{\bigcup\limits _{j\in \IN{p}}} \, K_j \py 
\IN{k} = \stackrel{_{_D}}{\bigcup\limits _{k\in \IN{p}}} \, J_k  \ .\earr
$$
\een
We remark that, if $\cF_0=\vacio$, 
 these facts are still valid for local minima 
by Theorem \ref{vale la conjetur}. 
The next three Propositions give a complete proof of Theorem \ref{remando}. 
The first of them justifies the convention that $\la = \la\ua(S_{\cF_0})$. 
\EOE
\end{num}

\begin{pro}\label{los la}
\ID \ with $\la = \la\ua(\,S_{\cF_0}\,)$, and consider the notations of 
\ref{notaciones conj loc}. 
If $p>1$, then 
$$
\barr{rl}
i \in K_1 &\implies 
i< j  \ ( \implies 
\la_i \le \la_j \,)  \peso{for every} 
j \in \bigcup\limits_{r>1} K_r = \IN{s_\cF} \setminus K_1 \ . \earr
$$ 
Inductively, by means of Remark \ref{induc}, 
we deduce that all sets $K_j$ consist on consecutive indexes, and that 
$K_i<K_j$ (in terms of their elements) if $i<j$. 

\end{pro}
\proof
Suppose that there are $i \in K_1$ and $j \in K_r$ 
(for some $r>1$) such that  $j<i$.  
Then $\la_j \le\la_i\,$ and $\mu _i \le \mu_j\,$. For $t>0$ very small, 
let $\mu_i(t) = \mu_i-t >0$ and 
$\mu_j(t) = \mu_j+ t $. Consider the vector $\mu(t) $  obtained by changing in $\mu$ 
the entries $\mu_i $ by $\mu_i(t)$ and $\mu_j $ by $\mu_j(t)$.  Observe that 
not necessarily $\mu(t)  = \mu(t)  \da$, but we are indeed 
sure that $c_1>c_r\,$. 

\pausa
Nevertheless, by Remark \ref{rem doublystohasctic},  
$(\mu _i \coma \mu_j) \prec 
(\mu _i(t) \coma \mu _j(t)\,)\implies \ca \prec \mu \prec \mu(t)$. Therefore there exists 
$\cF'= (\cF_0\coma \cG')\in \cC_\ca(\cF_0)$ such that, using the 
ONB of Eq. \eqref{la BON2bis},  
$$
\barr{rl} 
S_{\cG' } &= \suml_{h\in \IN{d}} \, \mu_h(t) \cdot v_h \otimes v_h 
\py 
S_{\cF'}=S_{\cF_0} + S_{\cG' } = \suml_{h\in \IN{d}} \, (\la_h + \mu_h(t)\,)
 \, v_h \otimes v_h \ . \earr
$$
Denote by $V=\gen\{v_i\coma v_j\}$, which reduces both $S_\cF$ and 
$S_{\cF'}\,$. Also $S_{\cF'}|_{V\orto} = S_{\cF}|_{V\orto}\,$.
Considering the restrictions to $V$ as operators in $L(V)
\cong \matrec{2}$ we get that 
$$
\la(S_{\cF'}|_{V}) = 
(\la_i + \mu_i(t) \coma \la_j + \mu_j(t)\,) = 
(c_1-t \coma c_r+t ) \prec (c_1\coma c_r) = \la(S_{\cF}|_{V})
\quad \mbox{strictly ,}
$$
for $t$ small enough in such a way that 
$c_1-t>c_r+t$, so that $(c_1-t \coma c_r+t ) = (c_1-t \coma c_r+t )\da$.
Then the map $F$ of Eq. \eqref{La F}, considered 
both on $\R_{\geq 0}^2$ and $\R_{\geq 0}^d\,$,  satisfies that
$$
F\big(\, \la(S_{\cF'}|_{V})\, \big) <
F\big(\, \la(S_{\cF}|_{V})\, \big) \implies 
P_f(\cF') 
= F\big( \,\la ( S_{\cF'} )\,\big) < F(\la(S_\cF)\,)= P_f(\cF) \ , 
$$
a contradiction. 
The inductive argument follows from Remark \ref{induc}. \QED

\begin{num}\label{local global}
In the following two statements we assume that, 
for some $f \in \convfs$,  
the sequence $\cF= (\cF_0\coma \cG)\in \cC_\ca^{\rm op}(\cF_0)$ is 
a {global} minimum for $P_f\,$, or it is a {local} 
minimum if $S_{\cF_0}=0$ and $\la = 0$.  In both cases 
\ref{notaciones conj loc} applies. 
\EOE
\end{num}

\begin{pro}\label{los J ordenados}
\ID , and let $\cF= (\cF_0\coma \cG)\in \cC_\ca^{\rm op}(\cF_0)$ 
as in \ref{notaciones conj loc}  and \ref{local global}. 
Suppose that $p>1$. 
Given $h \in J_i$ and $l \in J_r\,$ then 
$$
i<r \implies a_h-a_l \ge c_i - c_r > 0 \ .
$$
In particular, the sets $J_i$ consist of consecutive indexes, 
and $J_1<J_2<\ldots<J_p\,$ (in terms of their elements).  
\QED
\end{pro}
\proof
Let us assume that  $i<r  \in \IN{p}\,$, $h \in J_i$ and $l \in J_r\,$,  
but  $l<h$ (even less: 
that $a_l \ge a_h\,$). Then  
$$ 
g_l\otimes g_l \le S_\cG\le S_{\cF} \py S_\cF \ g_l= c_r\, g_l \implies 
a_h = \|g_h\|^2 \le\|g_l\|^2 = a_l \le  c_r <  c_i \ .
$$ 
We also know that $\langle g_l\coma g_h\rangle=0$. 
Denote by $w_h =\frac{g_h}{\|g_h\|}=a_h^{-1/2}\, g_h $ 
and  $w_l =\frac{g_l}{\|g_l\| }=a_l^{-1/2}\, g_l\ $. Let 
$$
g_h(t)=\cos(t)\ g_h + \sin(t)\ \|g_h\| \ w_l 
\py 
g_l(t)= \cos(\gamma t)\ g_l + \sin(\gamma t)\ \|g_l\| \ w_h 
\peso{for} 
 t\in \R 
$$ 
for some convenient $\gamma>0$ that we shall find later.
Let  $\cF_\ga(t)$ be the sequence obtained by changing in $\cF$ 
the vectors $g_h$ by $g_h(t)$ and $g_l$  by $g_l(t)$, for every  $t\in\R$. 
Notice that $\|g_h(t)\|^2= a_h$ and  $\|g_l(t)\|^2=a_l$ for every $t\in\R$, so that 
all the sequences $\cF_\ga(t) \in \cC_\ca(\cF_0)$. 

\pausa 
Let $W = \gen \{w_h \coma w_l\}$, a subspace which reduces 
$S_\cF$ and $S_{\cF_\ga(t)}\,$. Note that $g_h(t),\, g_l(t)\in W$. 
In the matrix representation  with respect to this basis of $W$ we get that 
$$ 
g_h\otimes g_h= \bm{cc}a_h &0 \\0&0\em \barr{l} w_h \\ w_l \earr \ , \ \ \ 
g_h(t)\otimes g_h(t)= a_h \ \ \bm{cc} \cos^2(t)  & \cos(t)\,\sin(t) \\ 
\cos(t)\,\sin(t)  &  \sin^2(t) \em  \barr{l} w_h \\ w_l \earr  \ ,
$$
$$
g_l\otimes g_l= \bm{cc}0 &0 \\0&a_l \em \barr{l} w_h \\ w_l \earr  \py 
g_l(t)\otimes g_l(t)= a_l \ \ \bm{cc} \sin ^2(\ga \,t)  & \cos(t)\,\sin(t) \\ 
\cos(t)\,\sin(t)  &  \cos^2(\ga \, t) \em  \barr{l} w_h \\ w_l \earr 
$$
If we denote by $S(t) = S_{\cF_\ga(t)}\,$,  we get that 
$$
S(t) = S_\cF - g_h \otimes g_h - g_l\otimes g_l + 
g_h(t)\otimes g_h(t) + g_l(t)\otimes g_l(t) 
\ .
$$
Therefore $S(t)|_{W^\perp} =S_\cF|_{W^\perp} \,$. 
On the other hand, $S_\cF|_{W} = \bm{cc}c_i&0\\0&c_r\em$. Then
$$
S(t)|_W =  
\bm {cc} c_i+ a_h \,(\cos^2(t)-1)+a_l\,\sin^2(\gamma t) & 
a_h \,\cos(t) \,\sin(t)+ a_l\,\cos(\gamma t)\,\sin(\gamma t) \\ 
a_h \,\cos(t) \,\sin(t)+ a_l\,\cos(\gamma t)\,\sin(\gamma t) & 
c_r + a_h  \,\sin^2(t) +a_l ^2\,(\cos^2(\gamma t)-1) \em \igdef A_\ga(t) \ .
$$
Note that 
$\tr \, A_\ga(t) = c_i + c_r$ for every $t\in \R$. 
Therefore $\la( A_\ga(t)\,) \prec 
(c_i \coma c_r)$ strictly $\iff \|A_\ga(t)\|_{_2} ^2 < c_i^2 + c_r^2\,$.
Hence we consider the map $m_\ga : \R \to \R $ given by 
$$
m_\ga(t) =    \|A_\ga(t)\|_{_2} ^2 = \tr \, (A_\ga(t)^2) 
\peso{for every} t\in \R  \ . 
$$
Note that 
$S(0) = S_\cF \implies m_\ga(0) = c_i^2 + c_r^2\,$. We shall see that, for a convenient choice 
of $\ga$, it holds that $m_\ga'(0) = 0$ but $m_\ga''(0) < 0$. This will contradict
the (local) minimality of $\cF$, because $m_\ga$ would have in this case 
a maximum at $t= 0$, so that $\la( A_\ga(t)\,) \prec 
(c_i \coma c_r)$ strictly 
$\stackrel{\eqref{cachos}}{\implies} \la (S_{\cF_\ga(t)}) \prec \la(S_\cF)$ strictly  $\implies P_f(\cF_\ga(t)\,) < P_f(\cF)$ for every $t$ near $0$. 

\pausa
Indeed, we first compute the derivatives of the entries $a_{ij}$ of $A_\ga(t)$ : 
$$
\barr{lcr}
a_{11}'&=& - a_h \,\sin(2t)+\gamma \ a_l\,\sin(2 \gamma t) \\
a_{12}' &=&a_h \,\cos(2t) + \gamma \ a_l\,\cos(2 \gamma t)\\
a_{22}'&=&a_h \,\sin(2t) - \gamma \ a_l\,\sin(2\gamma t)
\earr 
\py 
\barr{lcr} a_{11}''&=& 2\,[ - a_h \,\cos(2t)+\gamma^2 \ a_l\,\cos(2 \gamma t)] \\
a_{22}''&=&2\,[a_h \,\cos(2t) - \gamma^2 a_l\,\cos(2\gamma t)]
\earr \ . 
$$
So $ a_{11}'(0)=0\, ,\ a_{22}'(0)=0 $ and $ a_{12}(0)=0$.
Then, for $ i\coma j\in\IN{2}$ we have that 
$$
(a_{ij}^2)'(0)=2 \ a_{ij}(0)\ a_{ij}'(0)=0 \py 
(a_{ij}^2)''(0) = 2\big(\, (a_{ij}')^2(0)
+a_{ij}(0) \ a_{ij}''(0)\, \big) \ .
$$
Therefore 
$
(a_{11}^2)''(0) =  4c_i( - a_h+ \gamma^2 \ a_l)$, 
$(a_{12}^2)''(0) = 2 \ (a_h  + \gamma \ a_l)^2$ and 
$(a_{22}^2)''(0) = - 4 c_r \ (- a_h  + \gamma^2 \, a_l)\,$. 
We conclude that  $m_\ga'(0)=0$ (for every $\ga\in\R$) and that 
$$ 
m_\ga''(0)= 4\ \Big[\, c_i( - a_h+ \gamma^2 \ a_l) + (a_h  
+ \gamma \ a_l)^2 - c_r \ (- a_h  + \gamma^2 \, a_l)\, \Big] \ ,
$$ 
which is quadratic polynomial on $\ga$ with discriminant 
(if we drop the factor $4$) given by 
$$ 
D= a_h\ a_l\ \Big[\ a_h\ a_l - \big(\, 
a_l+(c_i-c_r)) \ (a_h- (c_i-c_r)\, \big) \ \Big]  \ .
$$ 
As we are assuming that $a_l \ge a_h\,$ then  $D>0$, because  
$$ 
\big(\, 
a_l+(c_i-c_r)) \ (a_h- (c_i-c_r)\, \big) 
= a_l\ a_h - (c_i-c_r) (a_l-a_h) - (c_i-c_r)^2
< a_l\ a_h \ .
$$
Hence there exists $\ga\in \R$ such that $m_\ga''(0) <0$. 
Observe that as long as $0<(c_i-c_r) (a_l-a_h) + (c_i-c_r)^2 \, (\iff 
a_h-a_l < c_i-c_r)\,$ we arrive at the same contradiction. 
\QED

\pausa
The following result is inspired on some ideas from \cite {casminpot}.

\begin{pro}\label{primeros LI}
\ID , and let $\cF= (\cF_0\coma \cG)\in \cC_\ca^{\rm op}(\cF_0)$ 
as in \ref{notaciones conj loc}  and \ref{local global}. 
For every  $ j<p$, the subsequence  
$\{g_i\}_{i\in J_j}$  of $\cG$ is linearly independent. 
\end{pro}
\proof
Suppose that there exists $j\in\IN{p-1}$ such that $\{g_i\}_{i\in J_j}$ is linearly dependent. 
Hence there exists coefficients $z_l\in \C$, $l\in J_j$ (not all 
zero) such that $|z_l|\leq 1/2$ and  
\begin{equation}\label{eq1}
\sum_{l\in J_j}\overline{z_l} \ a_l \ g_l=0\ . 
\end{equation}
Let $I_j \inc J_j$ be given by $I_j =\{l\in J_j:\ z_l\neq 0\}$
 and let $h\in \C^d$ such that $\|h\|=1$ and  $S_\cF h=c_p \,h$. 
For $t\in (-1,1)$ let $\cF(t)=(\cF_0,\cG(t))$ where  
$\cG(t)=\{g_i(t)\}_{i\in\IN{k}}$ is given by 
$$
g_l(t) = \begin{cases}  \ (1-t^2\,|z_l|^2)^{1/2} g_l+t\,z_l\,a_l h 
& \mbox{if} \ \ l\in I_j  \\
\quad \quad\quad g_l & \mbox{if} \ \ l\in \IN{k} \setminus I_j  \ \ .
\end{cases}  
$$
Fix $l\in I_j \,$. Let $\Preal(A)= \frac {A+A^*}{2}$ denote the real part of 
each  $A \in \op$. Then 
$$
g_l(t)\otimes g_l(t)
=(1-t^2\,|z_l|^2)\ g_l\otimes g_l+ t^2\,|z_l|^2\,a_l^2 \ h\otimes h 
+ 2\,(1-t^2\,|z_l|^2)^{1/2}\,t \ \Preal(h\otimes a_l\,z_l\, g_l)
$$
Let $S(t)$ denote the frame operator of $\cF(t)$ and notice that $S(0)=S_\cF$. 
Note that 
$$
S(t)=S_\cF+t^2 \sum_{l\in I_j }  |z_l|^2 \left( - g_l\otimes g_l 
+ a_l^2 \ h\otimes h \right) + R(t)$$ where 
$R(t)=2 \suml_{l\in I_j }(1-t^2\,|z_l|^2)^{1/2}\,t \ \Preal(h\otimes a_l\,z_l\, g_l)$. 
Then $R(t)$ is a smooth function such that 
$$
R(0) = 0 \ \ , \ \ 
R'(0)=\sum_{l\in I_j } \Preal(h\otimes a_l\,z_l\, g_l)
=\Preal(h\otimes \sum_{l\in I_j } a_l\,z_l\, g_l)=0\ ,
$$ 
and such that $R''(0)=0$. Therefore 
$\lim\limits_{t\rightarrow 0} \ t^{-2}\ R(t)=0 $.
We now consider 
$$
W=\gen\,\big(\,\{g_l:\ l\in I_j \}\cup \{h\}\,\big)
=\gen\,\big\{\,g_l:\ l\in I_j \,\big\}\perp \C\cdot h\ .
$$
Then $\dim W=s+1$,  for $s=\dim\gen\{g_l:\ l\in I_j \}\geq 1$. 
By construction, the subspace $W$ reduces $S_\cF$ and $S(t)$ for 
$t\in\R$, in such a way that $S(t)|_{W^\perp}=S_\cF|_{W^\perp}$ 
for $t\in \R$. On the other hand 
\beq\label{adet}
S(t)|_{W}=S_\cF|_W+t^2 \sum_{l\in I_j }  |z_l|^2 \left( - g_l\otimes g_l 
+ a_l^2 \ h\otimes h \right) + R(t) = A(t)+R(t)\in L(W)\ ,
\eeq
where we use the fact that the ranges of the selfadjoint operators 
in the second and third term in the formula above clearly lie in $W$. 
Then $\la\big(\,S_\cF|_W\,\big)=\big(\,c_j \, \uno_s\coma   c_p\, \big)
\in (\R^{s+1}_{>0})\da $ and 
$$\barr{rl}
\la \Big(\, \sum_{l\in I_j }  |z_l|^2  g_l\otimes g_l\,\Big) &
=(\gamma_1 \coma \ldots \coma \gamma_s \coma 0)
\in (\R^{s+1}_{\geq 0})\da \peso{with} \gamma_s>0 \ ,\earr
$$
where we have used the definition of $s$ and the fact that $|z_l|>0$ for $l\in I_j \,$. 
Hence, for sufficiently small $t$, 
the spectrum of the operator  $A(t)\in L(W)$
defined in \eqref{adet} is 
$$\barr{rl}
\la\big(\, A(t)\,\big) 
&=\big(\, c_j-t^2\,\gamma_s \coma \ldots \coma c_j-t^2 \,\gamma_1 \coma c_p
+t^2 \, \sum_{l\in I_j }a_l^2\,|z_l|^2 \,\big) \in (\R^{s+1}_{\geq 0})\da \ , \earr	
$$ 
where we have used the fact that $\langle g_l \coma h\rangle=0$ for every $l\in I_j \,$. 
Let us now consider 
$$
\la\big(\, R(t)\,\big)
=\big(\,\delta_1(t) \coma \ldots \coma \delta_{s+1}(t)\, \big) 
\in (\R^{s+1}_{\geq 0})\da\peso{for} t\in \R \ .
$$ 
Recall that in this case $\lim\limits_{t\rightarrow 0}t^{-2} \delta_j(t)=0$ 
for $1\leq j\leq s+1$. Using Weyl's inequality 
on Eq. \eqref{adet},  we now see that 
$\lambda \big(\,S(t)|_W\,\big)\prec \la\big(\, A(t)\,\big)
+\la\big(\, R(t)\,\big)\igdef \rho(t)\in (\R^{s+1}_{\geq 0})\da$. We know that 
$$
\barr{rl}
\rho(t)&= \big(\, c_j-t^2\,\gamma_s+\delta_1(t) \coma 
\ldots \coma c_j-t^2 \,\gamma_1+\delta_s(t) \coma 
c_p+t^2 \, \sum_{l\in I_j }a_l^2\,|z_l|^2  +\delta_{s+1}(t)\, \big) 
\\ &\\
&=  
\Big(\,c_j-t^2\,(\gamma_s-\frac{\delta_1(t)}{t^2}) \coma \ldots \coma 
c_j-t^2 \,(\gamma_1-\frac{\delta_s(t)}{t^2}) \coma 
c_p+t^2 \, (\sum_{l\in I_j }a_l^2\,|z_l|^2+\frac{\delta_{s+1}(t)}{t^2})\,\Big)  \ .
\earr
$$
A direct test shows that, for small $t$,  this $\rho(t) \prec 
\lambda(S_\cF|_W)=\big(\,c_j \, \uno_s\coma   c_p\, \big)$ strictly. 
Then, since $f$ is strictly convex, 
for every sufficiently small $t$ we have that  
$$
P_f\big(\,\cF(t)\,\big)\leq 
\tr \, f\big(\,\la(S_\cF|_{W^\perp})\,\big) +\tr \,f \big(\,\rho(t)\,\big)
< \tr \, f\big(\,\la(S_\cF|_{W^\perp})\,\big) +\tr \,f \big(\,\la(\,S_\cF|_{W}\,)\,\big) 
= P_f(\cF)  \ .
$$
This last fact contradicts 
the assumption that $\cF$ is a local minimizer of $P_f$ in $\cC_\ca^{\rm op}(\cF_0)$.
\QED

\begin{rem}\label{vieja conj}
Proposition \ref{primeros LI} allows to show that in case $\cF_0=\emptyset$ then local and global minimizers of a convex potential $P_f$, induced by $f\in \convfs$, on $\cC_\ca(\cF_0)$ - endowed with the product topology - coincide, as conjectured in \cite{MR}.


\pausa
Recall that a local minimizer  $\cF$ is a juxtaposition of tight frame sequences $\{\cF_i\}_{i\in \IN{p}}$ which generate pairwise orthogonal subspaces of $\hil$. Notice that by 
\cite[Lemma 4.9]{MRS5} $\cF$ is a frame for $\hil$. Moreover, by Proposition \ref{los J ordenados}, it is constructed using a partition of $\ca$ with consecutive indexes. 

\pausa
Now  by inspection of the proof of  Proposition \ref{primeros LI} we see that only one of such frame sequences  can be a linearly dependent set: that with the smallest tight constant $c_p$. This forces that the (ordered) spectrum $\nu$ of a local minimizer must be either  $\nu=c\uno_d$ or
$$
\nu=(a_1 \coma  a_2 \coma  \ldots \coma  a_r \coma c \coma \cdots \coma c)
 \ , \peso{where}
a_r>c\geq a_{r+1} \ , 
$$ 
and $c$ is the constant of the unique tight subframe constructed with a linear dependent sequence of vectors with norms given by $\{a_i\}_{i=r+1}^k$ (notice that this forces $c\geq a_{r+1})$. But it is not difficult to see that this vector can be constructed in a unique way, that is, there is only one $r$ such that 
$$
a_{r+1}\leq c=\frac{1}{d}\, \Big(\, \tr(\ca)-
\sum_{i=1}^ra_i\, \Big)<a_r \ .
$$ 
That is, the spectrum of local minimizers is unique and therefore local and global minimizers of $P_f$ coincide, for every potential $P_f$ as above.
\EOE
\end{rem}

\subsection{Several proofs.} \label{sect sev proof}
\ID \ with $\la = \la\ua(S_{\cF_0})$, $\ca = \ca\da$ and $d\le k$. 
Recall that we denote by  
$h_i=\la_{i}+a_i$ for every $i \in \IN{d}\,$ and, 
given  $j\coma r \in \IN{d}$ such that $ j\le r$, we denote by 
$$\barr{rl}
P_{j\coma r} &=\frac{1}{r-j+1}\ \suml_{i=j}^r\  h_i = 
\frac{1}{r-j+1}\ \suml_{i=j}^r\ \la_{i}+a_i 
\ .\earr
$$
We shall abbreviate $ P_{1\coma r}= P_r\,$.

\pausa
\begin{num}[{\bf Proof of Theorem \ref{remando}}]  \label{proof remando 2}
We rewrite its statement: 
\ID \ with $d\ge k$.  Let $\cF= (\cF_0\coma \cG)\in \cC_\ca^{\rm op}(\cF_0)$ be a {\bf global} 
minimum for $P_f\,$
on $\cC_\ca^{\rm op}(\cF_0)$.  
Using the notations of \ref{notaciones conj}, 
assume that $\la = \la\ua(S_{\cF_0})$, $\mu = \mu\da = \elmu$ and $\ca = \ca\da$. Then 
\ben
\item 
There exist indexes $0 = s_0<s_1 <\cdots< s_{p-1}<s_p 
= s_\cF= \max \{j \in \IN{d} : \mu_j \neq 0\}$ such that 
\beq\label{los Kj}
\barr{rl}
K_j &=J_j =\{s_{j-1}+1\coma \ldots \coma  s_j\} \ ,\quad 
\ ,\peso{for} j\in\IN{p-1}\  , \\&\\
K_p& =\{s_{p-1}+1 \coma \ldots \coma  s_p\} \ , \ J_p=\{s_{p-1}+1\coma \ldots\coma  k\} \ .
\earr 
\eeq
\item The vector $\nu_f(\la\coma \ca) = \big(\, c_1 \, \uno_{s_1} \coma 
\dots \coma c_p\, \uno_{s_p-s_{p-1}} \coma \la_{s_p+1} \coma 
\dots \coma \la_d\,\big)$, where 
\beq\label{los cr}
c_r=\frac{1}{s_r-s_{r-1}} \ \sum_{i=s_{r-1}+1}^{s_r}\, h_i
= P_{s_{r-1}+1\coma s_{r}}\peso{for} r\in\IN{p-1} \ ,
\eeq
or also $c_r = \la_j + \mu_j$ for every $j \in K_r=J_r $ for $ r\in\IN{p-1}\,$. 
\item The constant $c_p$ is the one defined by the feasible final part i.e.,  
 $c_p = Q_{s_{p-1}\coma  s_p}$ of \eqref{los Q bis} and the indexes $s_{p-1}$ and $s_p$ are determined by 
the last block, which is feasible. 
\een
\proof
Recall from Eq. \eqref{cajas} that for every $j \in \IN{p-1}$ 
$$
W_j \igdef \ker \,(S-c_j\,I_W\,) =\gen\{v_i : i \in K_j\} = \gen\{g_i : i \in J_j\} 
$$
By Proposition \ref {primeros LI} $|J_j| = \dim W_j = |K_j|$ for $j<p$.    
Using now Propositions \ref{los la} and \ref{los J ordenados},  
we deduce that there exist indexes $0 = s_0<s_1<  s_2<\cdots< s_{p-1}<s_p 
= s_\cF=\max \{j \in \IN{d} : \mu_j \neq 0\}$
such that the sets $K_j$ and $J_j$ satisfy Eq. \eqref{los Kj}. 
Using Eq. \eqref{cajas} again, 
\beq\label{mus as}
\barr{rl}
S_\cG|_{W_j} & = \sum_{i\in J_j} \, g_i  \otimes g_i 
\implies 
\tr \, S_\cG|_{W_j} =\sum_{i\in K_j} \, \mu_i = \sum_{i\in J_j} \, a_i   \ . \earr
\eeq
Therefore  $(s_j-s_{j-1}) \, c_j = \tr S|_{W_j} = 
 \tr \, S_{\cF_0}|_{W_i} + \tr \, S_\cG|_{W_i} =\sum_{i\in K_j} h_i \,$,
for every $j<p$. 
Then the vector $\nu_f(\la\coma \ca) = \big(\, c_1 \, \uno_{s_1} \coma 
\dots \coma c_p\, \uno_{s_p-s_{p-1}} \coma \la_{s_p+1} \coma 
\dots \coma \la_d\,\big)$, where 
the constants $c_r$ are given by Eq. \eqref{los cr} for $r <p$. 
Item 3 follows from Remark \ref{induc} and Lemma \ref{solo c1}. 
\QED
\end{num}

\begin{lem}\label{crece mayo}
\ID \ with $k\ge d$. Given $m\in \IN{d}\,$, 
$$
(a_j)_{j\in \IN{m}}   \prec 
(P_{ m}-\la_{j})_{j\in \IN{m}}
\iff P_m \geq P_i   \peso{for every} i\in \IN{m}
\iff P_{1\coma m} 
=\max\limits_{i\in \IN{m}}\{P_{1\coma i}\} \ .
$$ 
\end{lem}
\proof Straightforward. \QED

\begin{rem} \label{rm ck mayo}
\ID \, with $k\ge d$ and  recall the description of a minimum 
$\elnu$ given in 
Theorem \ref{remando}.  
As in Lemma \ref{crece mayo}
(or by an inductive argument using Remark \ref{induc}) 
we can assure that for every 
$r\le p-1$, the constants  
\beq\label{ck mayo}
c_r = P_{s_{r-1}+1\coma s_{r}} \ge P_{s_{r-1}+1\coma j } 
\peso{for every} j \peso{such that} s_{r-1}+1\le j \le s_r \ .
\eeq
It uses that $(a_j)_{j=s_{r-1}+1}^{s_r} \prec (\mu_j)_{j=s_{r-1}+1}^{s_r}
= (c_r-\la_{j})_{j=s_{r-1}+1}^{s_r}
\,$, a consequence of Eq. \eqref{mus as}. \EOE
\end{rem}

\begin{lem}\label{todo junto} 
\ID . With the notations of Theorem \ref{remando}, the global minimum $\elnu$, its  
constants $c_j$ and the indexes $s_j$ (for $j \in \IN{p}$) satisfy
the following properties: 
\ben
\item Suppose that $p>1$. For every $ j\in \IN{p-1}\,$
such that $j>1$,  the constant $c_j$ satisfies that 
\beq\label{ck y el pik}
c_j = P_{s_{j-1}+1\coma s_{j}} = 
\frac{1}{s_j-s_{j-1}}\sum_{i=s_{j-1}+1}^{s_j}\, h_i
< \frac{1}{s_j}\sum_{i=1}^{s_j}\, h_i = P_{1\coma s_j}  \ .
\eeq 
\item Fix $ j\in \IN{p-1}\,$ such that $j>1$. 
Then
\beq\label{P1t}
P_{1\coma t}<  P_{1\coma s_{j-1}}  \peso{for every} s_{j-1}<t\le s_{p-1} \ .
\eeq
\item In particular the averages 
$\ds P_{1\coma s_j} = \frac{1}{s_j}\sum_{i=1}^{s_j}h_i
<\frac{1}{s_{j-1}}\sum_{i=1}^{s_{j-1}}h_i = P_{1\coma s_{j-1}}\,$ 
for  $2\leq j\leq p-1$. 
\een
\end{lem}
\proof
The inequality of item 1 follows since 
$$\barr{rl}
\sum_{i=1}^{s_j}h_i&=\sum_{i=1}^{s_1}h_i
+\sum_{i=s_1+1}^{s_2}h_i+\ldots+\sum_{i=s_{j-1}+1}^{s_j}h_i\\&\\
&=s_1\, c_1+(s_2-s_1)\, c_2+\ldots +(s_j-s_{j-1})\,c_j>s_j\,c_j \ .
\earr$$
Now we prove the inequality of Eq. \eqref{P1t}: Given an index $t$ such that $s_{j-1}<t\le s_j\,$, 
$$
\barr{rl}
t\, P_{1\coma t} & = s_{j-1} \, P_{1\coma s_{j-1}} + 
\sum_{i=s_{j-1}+1}^t h_i \\&\\& = s_{j-1}\, 
P_{1\coma s_{j-1}} + (t-s_{j-1}) \ \frac{1}{(t-s_{j-1})} \ 
 \sum_{i=s_{j-1}+1}^t h_i \\&\\
&\stackrel{\eqref{ck mayo}}{\le}  s_{j-1} \, P_{1\coma s_{j-1}} +(t-s_{j-1})\, c_j \\&\\
&< s_{j-1}\,  P_{1\coma s_{j-1}} +(t-s_{j-1})\, c_{j-1} \\&\\
&\leq s_{j-1}\,  P_{1\coma s_{j-1}}+(t-s_{j-1}) \, P_{1\coma s_{j-1}} 
= t \, P_{1\coma s_{j-1}} \ , \earr
$$ 
where we used the fact that $c_{j-1}\, \leq P_{1\coma s_{j-1}} $ for $1\leq j-1\leq p-1$, 
which follows from item 1.  In particular we have proved item 3, and this also proves 
that Eq. \eqref{P1t} holds for $s_{j}<t\le s_{p-1}\,$.
\QED

\begin{pro} \label{el s_1}
\ID . With the notations of Theorem \ref{remando}, the global 
minimum $\nu = \elnu$, its  
constants $c_j$ and the indexes $s_j$ (for $j \in \IN{p}$) satisfy
the following properties:  Suppose we know the index $ s_{p-1}\,$, and that $p>1$. 
Then we have a recursive method to reconstruct $\nu$: 
\ben
\item The index $s_1 = \max \, \big\{j \le s_{p-1} \, :\, 
P_{1\coma j} = \max\limits_{i\le s_{p-1}}  \, P_{1\coma i} \, \big\}$, and 
$c_1 = P_{1\coma s_1}\,$.
\item If we already compute the index $s_j$ and $s_j<s_{p-1}\,$, then
$$
s_{j+1} = \max \, \big\{s_j< r \le s_{p-1} \, :\, 
P_{s_j+1\coma r} = \max\limits_{s_j< i\le s_{p-1}}  \, P_{s_j+1 \coma i} \, \big\} 
\py c_{j+1} = P_{s_j+1\coma s_{j+1}}\ .
$$
\een
\end{pro}
\proof
The formula $P_{1\coma s_1} = \max\limits_{i\le s_{p-1}}  \, P_{1\coma i}$
follows from Lemma \ref{crece mayo} and Eq. \eqref{P1t} 
of Lemma \ref{todo junto}, which also implies that 
$s_1$ must be the greater index (before $s_{p-1}$) satisfying this property. 

\pausa
The iterative program works by applying the last fact to 
the successive truncations of $\nu$ which are still minima 
in their neighborhood, by Remark \ref{induc}. \QED

\pausa
Recall that $h_i = \la_i + a_i$ and that, 
for $0\leq j< r\leq d$, we denoted by
\beq\label{los Q}
\barr{rl}
Q_{j\coma r} &=\frac{1}{r-j} \ \Big[\  \sum_{i=j+1}^{r}\  h_i 
+ \sum_{i=r+1}^k\  a_i  \Big] 
=\frac{1}{r-j} \ \Big[\  \sum_{i=j+1}^{k}\  a_i  + 
\sum_{i=j+1}^{r}\  \la_i    \Big] \ ,\earr
\eeq
and we abbreviate $Q_{1\coma r} =Q_r\,$. Recall also the notion of feasible indexes 
given in Definition \ref{s feas}:
\ID \ with $\la = \la\ua(S_{\cF_0})$ and $k\ge d$.
 Given $s\in \IN{d-1}$ denote by 
$\la ^s = (\la_{s+1} \coma \dots \coma 
\la_d )\in \R^{d-s}$ and $ \ca^s = (a_{s+1}\coma \dots \coma a_k)
$, the truncations of the original vectors $\la$ and $\ca$. 
Recall that  the index $s$ is  feasible if the pair 
$(\la^s\coma \ca^s)$ is feasible for the CP. 
In any case we denote by 
$$
\nu_s= \nu (\la^s\coma \ca^s) = 
\big(\,c\, \uno_{r-s}\coma   \la_{r+1}  \coma \dots \coma \la_d  \big)
\peso{where} c = Q_{s \coma r} 
$$
for the unique $r>s$ such that $\la_{r}\le c<\la_{r+1}\,$. 
This means that $\la_s\leqp \nu_s\in (\R^{d-s}_{_{>0}})\ua$ and that 
$\tr \nu_s = \tr\, \la^s + \tr \, \ca^s\,$.

\begin{lem} \label{el r de M}
\ID . Fix an index $s \in \IN{d-1}\cup \trivial \,$. Then 
\ben
\item
The index $r$ associated to $\nu_s$ as in 
the previous notations is given by 
$$
r = \max \ \{ w\in \IN{d}: w>s \py  Q_{s\coma w} = 
\min\limits_{j>s}\  Q_{s\coma j} \ \} \ .
$$
In other words, $r$ is the unique 
index which satisfies:  
Given $j >s$, 
\beq\label{el nuevo r}
Q_{s\coma r} <Q_{s\coma j} \peso{if} j>r \py
Q_{s\coma r} \le Q_{s\coma j} \peso{if} j< r \ .
\eeq
\item 
Given an index  $l\in \IN{d-1} \,$,  
\beq\label{menor que r}
l> s \py Q_{s\coma l} < \la_{\,l+1} \implies l\ge r \ ,
\eeq
where $r$ is the index associated to $\nu_s$ of item 1. 
\een
\end{lem}
\proof 
Item 1 follows from Proposition \ref{el r de M tris} 
applied to $\la^s$ and $\ca^s$. 

\pausa
Item 2 : Assume that $l <l+1\le  r$. 
Then $Q_{s\coma l}< \la_{\,l+1} \le \la_{r} \le Q_{s\coma r}\,$. 
In this case 
$$
\barr{rl}
\tr \, \la^s + \tr\, \ca^s & \stackrel{\eqref{los Q}}{=}  
(l-s)\, Q_{s\coma l}  \ 
+ \suml_{i=l+1}^d \la_i \\&\\ 
&=   (l-s)\, Q_{s\coma l}  \ + \suml_{l+1 \le i \le r} \la_i  
+ \suml_{i=r+1}^d \la_i \\&\\&
< (r-s)\, Q_{s\coma r} +  \suml_{i=r+1}^d \la_i
\stackrel{\eqref{los Q}}{=} \tr \, \la^s + \tr\, \ca^s
\ , 
\earr
$$
a contradiction. Hence $l\ge r $. \QED

\begin{pro}\label{el c es menor}
\ID \, which is not feasible, with $k\ge d$. Let 
$$
s^*= \min \ \{\ s\in \IN{d} \ : \  s  \ \mbox{ \rm is feasible } \} \ .
$$ 
Let $\nu^*$ be constructed using the recursive 
method of Proposition \ref{el s_1}, by using  $s^*$ instead of $s_{p-1}$ 
(which can always be done). Then if we get the constants 
$c_1>\ldots>c_{q-1}\,$, and we define $c_q$ as the feasibility
constant of $\lambda^{s^*}$ and  $\ca^{s^*}$, then  $c_{q-1}>c_q\,$.
\end{pro}

\proof For simplicity of the notations, 
by working with the pair $(\la^{s_{q-2}}\coma \ca^{s_{q-2}})$, 
 we can assume 
that $q =2$. Denote by $ s_1=s^*< s_2$ and $c_1\coma c_2$ 
the indexes and constants given by:
\beq\label{c1c2}
c_1=\frac{1}{s_1}\sum_{i=1}^{s_1} h_i=P_{1 \coma s_1} \peso{and} c_2
=Q_{s_1\coma s_2} 
=\frac{1}{s_2-s_1}\left(\sum_{i=s_1+1}^{s_2} h_i+\sum_{i=s_2+1}^k a_i\right) \ ,
\eeq
and we must show that $c_1>c_2\,$. 
Recall that $h_i=\la_{i}+a_i\,$. We can assume that:
\begin{itemize}
\item By Proposition \ref{el s_1},  $c_1\geq \frac{1}{p} \ \suml_{i=1}^{p} h_i=P_{1 \coma p}$ for 
every $p\in \IN{s_1}\,$.
\item $c_2\geq \frac{1}{p-s_1} \ \suml_{i=s_1+1}^{p} h_i=P_{s_1+1 \coma p}$ 
for every $s_1+1\leq p\leq s_2\,$.
\item $\la_{s_2}\leq c_2 <\la_{s_2+1}\ ,$
\end{itemize}
where the second item follows by the feasibility of $s^*$ 
and the last item says that $c_2$ is the feasible constant for the second block.

\pausa
Suppose that $c_1\leq c_2$ and we will arrive to a contradiction 
by showing that, in such case, the pair $(\la\coma \cb)$ 
would be  feasible (that is, $s^*=0$ or $s_{q-2}$).
In order to do that, let 
$$
\barr{rl}
t\in \IN{d} &\py 
b\igdef 
Q_t = \frac{1}{t} \ \left(\ \sum_{i=1}^t h_i +\sum_{i=t+1}^k a_i \ \right) 
\earr
$$ 
be the unique constant such that $ \la_{t}\leq b <\la_{t+1}\,$, which 
appears in $\nuel$. Then   
$$
\barr{rl}
c\igdef Q_{s_2} &= \frac{1}{s_2}\left(\suml_{i=1}^{s_2} h_i 
+\suml_{i=s_2+1}^k a_i\right) =
 \frac{1}{s_2}(s_1\, c_1+(s_2-s_1)\, c_2)
\leq c_2<\la_{s_2+1} \ .\earr
$$ 
By Eq. \eqref{menor que r} we can deduce that $t\leq s_2\,$. Moreover, 
by item 1 of Lemma \ref{el r de M} we know that 
\begin{equation}\label{feasible pierde}
\barr{rl}
b=Q_t &= \frac{1}{t}\left(\suml_{i=1}^t h_i +\suml_{i=t+1}^k a_i\right)
\leq \frac{1}{p}\left(\suml_{i=1}^p h_i +\suml_{i=p+1}^k a_i\right)
= Q_{p} \peso{ for every} p\in \IN{d} \ .\earr
\end{equation}
In particular, $b\leq c\leq c_2\,$. On the other side, $c_1\leq b$. Indeed, if $\nu = \nuel$ then 
$$
\la \leqp \nu^* \py t= \tr \, \nu^* =   \tr \, \nu 
\implies \nu \prec \nu^* \implies 
b  = \nu_1 \, \ge  \, \nu^*_1= c_1 \ ,
$$
because $c_1 \le c_2\implies \nu^* =(\nu^*)\ua$ 
and since $\nu = \nu\ua $ is the $\prec$-minimum of the set 
$$
\{\la\ua (S): S_{\cF_0} \le S \py \tr \, S = t\}   = 
\{\rho = \rho\ua : \la\leqp \rho \py \tr\,\rho=t\} \ , 
$$
by the remarks at the beginning of Section \ref{TFC} 
and Proposition \ref{el r de M tris}. 

\pausa
To show the feasibility, by Lemma \ref{crece mayo} 
we must show that 
$b\geq P_{1\coma p}$ for every $ p\in\IN{t}\,$.
First, if we are in the case $t\leq s_1\,$, this is clear since 
$b\geq c_1\geq P_{1 \coma p}$ for every $p\leq s_1\,$.
Finally, suppose that $t\geq s_1+1$. As before, 
$b\geq c_1$ implies $b\geq P_{1 \coma p}$ for every  $ p\leq s_1\,$. 
On the other hand, if $s_1< p\le t$ then 
Lemma \ref{el r de M} applied to $\nu_{s_1}$ (whose $``r"$ is $s_2\,$) assures that 
$$
\barr{rl}
c_2 < Q_{s_1\coma t} &\implies 
(t-s_1)\,c_2\leq \suml_{i=s_1+1}^t h_i +\suml_{i=t+1}^k a_i
\stackrel{\eqref{c1c2}}{=} t\, b-s_1\, c_1  \ . \earr
$$ 
Since $p\leq t$ and $b\le c_2\,$,  
this implies that $(p-s_1)\, c_2\le  p\,b-s_1\,c_1\,$. 
Therefore
\beq
p\,P_{1\coma p}= s_1\,c_1+(p-s_1)\, P_{s_1+1 \coma p} 
\leq s_1\,c_1+(p-s_1)\,c_2 \leq p\,b \ .
\QEDP
\eeq

\begin{pro}\label{s = sf}
\ID . With the notations of Theorem \ref{remando}, 
the global minimum $\elnu$ satisfies that 
$$
s_{p-1} = \min \ \{\ s\in \IN{d} \ : \  s  \ \mbox{ \rm is feasible } \} \ .
$$
\end{pro}
\proof
Denote by $s^*  $ the minimum of the statement. 
Since $s_{p-1}$ is feasible (recall the remark after 
Definition \ref{s feas}), then $s^*  \le s_{p-1}\,$. 
On the other hand, let us construct the  vector $\nu^*$ 
of Proposition \ref {el c es menor}, using 
the iterative method of Proposition \ref{el s_1} with respect 
to the index $s= s^* \,$, and  the solution 
for the feasible pair $(\la^{s^* }\coma \ca^{s^* })$ 
after $s^* \,$. Write $\nu^* = (\nu^*_1\coma \dots \coma \nu^*_s 
\coma c \, \uno_{r-s} \coma 
\la_{r+1} \coma \dots \coma \la_d)$, where $c$ is the 
 constant of the feasible part of $\nu^*$. Observe that 
Proposition \ref {el c es menor} assures that 
$c< \min \{\nu^*_i : 1\le i\le s\}$. 

\pausa
Using this fact and Proposition  \ref{el s_1}
it is easy to see that the vector $\mu = \nu^* - \la\ua $
satisfies that $\mu=\mu\da$. 
On the other hand Lemma \ref{crece mayo} and Remark \ref{rm ck mayo} 
assure that $\ca \prec \mu$ (using the majorization in each block 
and the fact that $\ca = \ca\da$). 
Then $\mu \in \Gamma_d(\ca )$ and $\nu^* \in \Lambda_\ca^{\rm op}(\lambda)$, 
the set defined in  Eq. \eqref{el convex}. 

\pausa
Moreover, in each step of the construction
of the minimum $\nu = \elnu$ we have to get the same index $s_j = 
s_j(\nu^*)$ of $\nu^*$ or there exists a step 
where the maximum which determines $s_j$   (for $\elnu$) 
satisfies that $s_j > s^* \,$ (in the eventual case 
in which $s_{p-1}>s^* $). 

\pausa
In both cases, we get that $\nu^*_i \le \nu_i$ for every 
index $1\le i \le s^* \,$. 
Consider the subvector of $\nu^* $ given by 
$\rho = (\nu^*_1\coma \dots \coma \nu^*_s \coma 
\la_{r+1} \coma \dots \coma \la_d) \in \R^{s+d-r}$, and
the respective part of $\elnu$ given by 
$\xi = (\nu_1 \coma \dots \coma \nu_s \coma \nu _{r+1} \coma \dots \nu_d)$. 
Since $\tr \, \nu^*  = \tr \,\nu$, the previous remarks show that 
$$
\rho \leqp \xi \implies \rho \prec_w \xi 
\implies (\rho \coma c \uno_{r-s} ) 
\prec (\xi\coma \nu_{s+1}\coma \dots \coma \nu_{r}) \  , 
$$
where the final majorization follows using  
Lemma 4.6 of \cite{MRS4}, which can be used since 
the constant $c< \min \{\nu^*_i : 1\le i\le s\}$ 
by Proposition \ref {el c es menor} (and because $c<\la_{r+1}$). 
Since majorization is invariant under rearrangements, we deduce that 
$\nu^* \prec \nu$. 

\pausa
Finally, using  Theorem \ref{teo sobre unico espectro}  
we know that $\nu= \elnu$ is the unique minimum for the map 
$\tr\, f(\cdot)$ in the set $\Lambda_\ca^{\rm op}(\lambda)$.  
This implies that 
$\nu^*=\nu$, and therefore $ s_{p-1} = s^* \,$. \QED

\fontsize {8}{9}\selectfont

\end{document}